%% file: main.tex
\documentclass[11pt,a4paper]{article}
\usepackage[letterpaper,left=1in,right=1in,top=1in,bottom=1in,includehead=true]{geometry}
\usepackage[utf8]{inputenc}
\usepackage[T1]{fontenc}
\usepackage{graphicx} 
\usepackage{caption}
\usepackage{subcaption}
\usepackage{amsmath}
\usepackage{amssymb}
\usepackage{amsthm}
\usepackage{mathtools}
\usepackage[ruled,vlined,linesnumbered]{algorithm2e}
\usepackage[dvipsnames]{xcolor}
\usepackage{subcaption}
\usepackage{authblk}
\usepackage{boldline}
\usepackage{tikz}
\usepackage{siunitx}
\usepackage{stmaryrd}
\usepackage{fancyhdr}
\usepackage{titlesec}
\usepackage{secdot}
\usepackage{multirow}
\usepackage{hyperref}
\usepackage{cleveref}
\usepackage{algpseudocode}
\usepackage{lineno}
\hypersetup{
    colorlinks=true,
    linkcolor=blue,
    citecolor=blue,
    filecolor=magenta,      
    urlcolor=cyan,
}
\usetikzlibrary{calc,angles,positioning,intersections,quotes,3d,shadings,arrows.meta}

\titleformat{\section}{\bfseries\large}{\thesection.}{0.5em}{}
\titleformat{\subsection}{\large}{\thesubsection.}{0.5em}{\itshape}

\newcommand{\Q}{\mathbf{U}}
\newcommand{\Sc}{\mathbf{S}}
\newcommand{\Fc}{\mathbf{F}}
\newcommand{\Gc}{\mathbf{G}}

\newcommand{\Fchat}{\widehat{\Fc}}
\newcommand{\Gchat}{\widehat{\Gc}}

\newcommand{\ten}[1]{\mathcal{#1}}
\newcommand{\mat}[1]{\mathbf{#1}}

\title{High-order Tensor-Train Finite Volume Method for Shallow Water Equations}
\author[a]{M. Engin Danis}
\author[a]{Duc P. Truong}
\author[b]{Derek DeSantis}
\author[b]{Mark Petersen}
\author[a]{Kim {\O.} Rasmussen}
\author[a]{Boian S. Alexandrov}
\affil[a]{Theoretical Division, Los Alamos National Laboratory, Los Alamos, NM 87545, USA}\affil[b]{Computer, Computational and Statistical Sciences Division, Los Alamos National Laboratory, Los Alamos, NM 87545, USA}
\date{\today}
\begin{document}
%
\maketitle
\begin{abstract}
    In this paper, we introduce a high-order tensor-train (TT) finite volume method for the Shallow Water Equations (SWEs). We present the implementation of the $3^{rd}$ order Upwind and the $5^{th}$ order Upwind and WENO reconstruction schemes in the TT format. It is shown in detail that the linear upwind schemes can be implemented by directly manipulating the TT cores while the WENO scheme requires the use of TT cross interpolation for the nonlinear reconstruction. In the development of numerical fluxes, we directly compute the flux for the linear SWEs without using TT rounding or cross interpolation. For the nonlinear SWEs where the TT reciprocal of the shallow water layer thickness is needed for fluxes, we develop an approximation algorithm using Taylor series to compute the TT reciprocal. The performance of the TT finite volume solver with linear and nonlinear reconstruction options is investigated under a physically relevant set of validation problems. In all test cases, the TT finite volume method maintains the formal high-order accuracy of the corresponding traditional finite volume method. In terms of speed, the TT solver achieves up to 124x acceleration of the traditional full-tensor scheme.
\end{abstract}
%
\input{text.tex}
\clearpage
\bibliographystyle{plain}
\bibliography{references}
\end{document}

%% file: text.tex
%
%
\section{Introduction}\label{sec:introduction}

As computer architectures evolve, new algorithms are often needed to make the best use of the latest equipment.  For example, in the 1990s there was a major transition from vector supercomputers to distributed memory clusters, and internode communication became the slowest part of a simulation. In global atmospheric models this spurred a transition from spectral methods\cite{kiehl_ccm3}, which require global communications for the transforms between physical and spectral space, to methods such as finite volume\cite{thuburn_numerical_2009} and spectral element\cite{dennis_high-resolution_2005}, which only require local communication with nearby processors. 

In recent years there has been another major change in the landscape of supercomputer architectures, as the main driver for sales became machine learning (ML) and artificial intelligence (AI) applications. 
Tensor cores are a new generation of chips specifically designed for AI and deep learning, such as the NVIDIA Volta, Turing, and Ampere classes of GPUs. 
Because AI, rather than computational physics, is driving the direction of commodity architecture, developers of numerical methods must seek out algorithms that are performant on this new hardware. The AI revolution has also spurred the development of libraries that are specifically tuned to run AI algorithms, such as Neural Networks, on tensor cores. Here we present recently developed numerical methods that leverage tensor networks. These methods manipulate large-scale data based on generalizations of the singular value decomposition to higher dimensional tensor arrays. 


A second driver of new algorithm development is to greatly increase the speed and resolution of global climate simulations. High-resolution global simulations are now typically 6km to 10km in the ocean and atmosphere \cite{caldwell_doe_2019} and the latest cloud-resolving simulations are at 3km resolution\cite{scream_2024}. This involves millions of horizontal cells and up to 128 vertical layers. In addition, climate research requires long simulations for spin-up, and ensembles of simulations to investigate the intrinsic variability of the climate system\cite{Kay_2015}. All these factors taken together produce the ``curse of dimensionality'', where simulation campaigns in climate science require many months on large supercomputers.

Tensor Networks (TNs)\cite{cichocki2014tensor} are a promising new approach that mitigate the effects of the curse of dimensionality, and may also take advantage of specialized AI hardware. TNs are a generalization of matrix factorizations to higher dimensions, whereby multidimensional data structures (tensors) are decomposed into manageable blocks. The computations normally performed on the large dataset (e.g. finite differences) can alternatively be performed on these smaller tensor components, drastically reducing computational costs. The most popular TN method is known as the Tensor Train (TT). In TT format, the large dataset is decomposed into a sequence of lower-dimensional tensors (carts), linked together in a chain (train), to efficiently represent and manipulate large-scale data \cite{oseledets2010tt}.  

TN techniques show great promise to attack the curse of dimensionality across a large range of problems.  Just in the last few years, tensor methods have been used to model the Navier-Stokes equations in a number of standard test cases: A backward-facing step \cite{demir2024tensor}, a lid-driven cavity \cite{kiffner2023tensor}, a turbulent flow in a channel \cite{vonLarcher2019tensor}, three-dimensional Taylor-Green vortex \cite{gourianov2022tensor}.  In recent work, we used TT methods to accelerate compressible flow simulations by up to 1000 times \cite{danis2024tensortrain}. These successes show that a tensor-based approach to fluids problems is both possible and promises greatly increased model efficiency. However, to date \emph{no work has been done to extend these methods to models of oceanic or atmospheric circulation.} Successful compression and speedup of geophysical fluid simulations via TN has the promise to radically alter the landscape of weather and climate modeling, allowing researchers to explore higher resolutions and larger ensembles.

Any new numerical method for atmosphere and ocean models must progress through a sequence of verification steps in order to be accepted by the community. The Shallow Water Equations (SWE) are a reduced equation set used as a first step for modeling geophysical fluids\cite{thuburn_numerical_2009,weller_voronoi_2009,Archilbald_2011,Lilly_2023}. They contain the relevant dynamics of atmospheric and oceanic flows: the Coriolis force and pressure gradient term for geostrophic balance, as well as horizontal advection of momentum and mass. At the same time, the SWEs are simple enough for rapid code development. Critically, the SWE may be tested against exact solutions during model development\cite{bishnu2024}, which is not the case for the following levels of complexity in the development sequence\cite{petersen_evaluation_2015}. The key assumptions of the SWE are that horizontal scales of motion are much larger than the vertical, which means it is hydrostatic, and that the fluid is incompressible so that it has uniform density\cite{cushman2011introduction}. This conveniently avoids the complexity of the equation of state, moist dynamics and clouds in the atmosphere, and vertical advection. Published test sets using the SWE have become an essential component of model development and verification\cite{williamson1992standard,calandrini_comparing_2021}. 

This article assesses the computational advantages of TN in modeling the SWE across a range of test cases.  In particular, we focus on TT methods for high-order finite volume methods for the SWE. The paper is organized as follows.  In Section \ref{sec:review-of-num-methods}, we review the SWE and our high-order methods for solving them.  Section \ref{sec:TT-decomposition} we cover the basics of the TT decomposition, and in Section \ref{sec:tensorization-of-the-method} we discuss how the finite volume scheme can be formulated in the TT format. We end with Section \ref{sec:results} covering the results for a series of test cases.

\section{Governing Equations and the Numerical Method}\label{sec:review-of-num-methods}
In this section, we will review the Shallow Water Equations (SWEs) and the finite volume method used to solve these equations. This discussion will only involve essential information for a typical finite volume implementation on traditional grids to lay the ground for the tensor-train implementation.
 
\subsection{Shallow Water Equations}\label{sec:governing-equations}
In this study, we consider both linear and nonlinear SWEs with a flat bottom topography. Both equations are solved in the conservative form,
\begin{equation}\label{eq:SWE-system}
    \frac{\partial \Q}{\partial t} + \frac{\partial \Fc}{\partial x} + \frac{\partial \Gc}{\partial y} = \Sc,
\end{equation}
where $\Q$ is the vector of conserved variables, $\Fc,\Gc$ are the fluxes in $x$ and $y$ directions, and $\Sc$ is the source term.

\emph{In the linear case}, \Cref{eq:SWE-system} is solved with 
\begin{equation}\label{eq:linear-SWE}
    \Q = 
    \begin{pmatrix}
        \eta \\
        u    \\
        v    \\
    \end{pmatrix},
    \qquad
    \Fc = 
    \begin{pmatrix}
        Hu \\
        g\eta  \\
        0      \\
    \end{pmatrix},
    \qquad
    \Gc = 
    \begin{pmatrix}
        Hv \\
        0      \\
        g\eta  \\
    \end{pmatrix},
    \qquad
    \Sc = 
    \begin{pmatrix}
        0     \\
        fv    \\
       -fu    \\
    \end{pmatrix},
\end{equation}
where the vector of conserved variables, $\Q$, consists of the surface elevation $\eta$, the $x$-velocity $u$ and the $y$-velocity $v$. Furthermore, $H$ is the mean depth of the fluid at rest, $f$ is the Coriolis parameter, and $g$ is the acceleration of gravity.

\emph{In the nonlinear case}, \Cref{eq:SWE-system} is solved with 
\begin{equation}\label{eq:nonlinear-SWE}
    \Q = 
    \begin{pmatrix}
        h \\
        hu    \\
        hv    \\
    \end{pmatrix},
    \qquad
    \Fc = 
    \begin{pmatrix}
        hu \\
        hu^2+\frac{1}{2}gh^2  \\
        huv      \\
    \end{pmatrix},
    \qquad
    \Gc = 
    \begin{pmatrix}
        hv \\
        huv      \\
        hv^2+\frac{1}{2}gh^2  \\
    \end{pmatrix},
    \qquad
    \Sc = 
    \begin{pmatrix}
        0   \\
        fhv \\
       -fhu \\
    \end{pmatrix},
\end{equation}
where $h$ is the shallow water layer thickness above the bottom topography.

%
%
%
%
\subsection{High-order Finite Volume Method for hyperbolic conservation laws}\label{sec:fv-review}
For simplicity, let us consider the 2-dimensional scalar hyperbolic conservation law,
\begin{equation}\label{eq:2d-scalar-h-pde}
    \frac{\partial u}{\partial t}+\frac{\partial f(u)}{\partial x}+\frac{\partial g(u)}{\partial y}=0,
\end{equation}
where $u$ is a generic conserved variable (not to be confused with the $x-$velocity), $f(u)$ and $g(u)$ are fluxes in $x$ and $y$ directions, respectively and we assume that proper initial and boundary conditions are provided. On a uniform mesh with grid spacing $\Delta x$ and $\Delta y$ in $x$ and $y$ directions, a finite volume method solves \Cref{eq:2d-scalar-h-pde} for the cell averages of $u$ in a given cell $(i,j)$,
\begin{equation}\label{eq:cell-avg-def}
    \widetilde{\overline{u}}_{i,j} = \frac{1}{\Delta x\Delta y}\int_{x_{i-1/2}}^{x_{i+1/2}}\int_{y_{j-1/2}}^{y_{j+1/2}}u\,dxdy,
\end{equation}
where $\overline{u}$ denotes cell average in $x$ and $\widetilde{u}$ denotes cell average in $y$. 

In this spirit, the semi-discrete cell-averaged form of \Cref{eq:2d-scalar-h-pde} for a cell $(i,j)$ is given as
\begin{equation}\label{eq:2d-scalar-h-pde-fv}
    \begin{aligned}
    \frac{d\widetilde{\overline{u}}_{i,j}}{dt}&+\frac{1}{\Delta x\Delta y}
    \int_{y_{j-1/2}}^{y_{j+1/2}} \left(f\left(u\left(x_{i+1/2},y\right)\right)-f\left(u\left(x_{i-1/2},y\right)\right)\right)\,dy \\
    &+\frac{1}{\Delta x\Delta y}
    \int_{x_{i-1/2}}^{x_{i+1/2}} \left(g\left(u\left(x,y_{j+1/2}\right)\right)-g\left(u\left(x,y_{j-1/2}\right)\right)\right)\,dx =0,
    \end{aligned}
\end{equation}
which can be approximated as
\begin{equation}\label{eq:2d-scalar-h-pde-fv-2}
    \frac{d\widetilde{\overline{u}}_{i,j}}{dt}+\frac{\widehat{f}_{i+1/2,j}-\widehat{f}_{i-1/2,j}}{\Delta x}+\frac{\widehat{g}_{i,j+1/2}-\widehat{g}_{i,j-1/2}}{\Delta y}=0,
\end{equation}
where the numerical fluxes $\widehat{f}_{i\pm1/2,j}$ and $\widehat{g}_{i,j\pm1/2}$ approximate the surface integrals by the 1-dimensional Gauss-Legendre quadrature rule with quadrature points $\delta_m$ and weights $w_m$:
\begin{equation}
    \begin{aligned}
        \widehat{f}_{i\pm1/2,j}&=\sum_{m}w_m \widehat{f}\left(u_{i\pm1/2,y_j+\delta_m\Delta y}^-,u_{i\pm1/2,y_j+\delta_m\Delta y}^+\right), \\
        \widehat{g}_{i,j\pm1/2}&=\sum_{m}w_m \widehat{g}\left(u_{x_i+\delta_m\Delta x,j\pm1/2}^-,u_{x_i+\delta_m\Delta x,j\pm1/2}^+\right),
    \end{aligned}
\end{equation}
where $\widehat{f}(u^-,u^+)$ and $\widehat{g}(u^-,u^+)$ denote the local Lax-Friedrichs flux. In the $x$-direction, for example, the local Lax-Friedrichs flux is defined as
\begin{equation}\label{eq:LLF}
    \widehat{f}(u^-,u^+) = \frac{1}{2}\left(f(u^-)+f(u^+)\right)-\frac{\lambda}{2}\left(u^+-u^-\right),
\end{equation}
where $u^\pm$ are point-wise values at the cell interfaces approximated by a high-order reconstruction method from the cell-averages $\widetilde{\overline{u}}_{i,j}$ and $\lambda=\max_{u\in\left(u^-,u^+\right)} |f'(u)|$.
%
%
%
%
\subsection{High-order Reconstruction}\label{sec:fv-recon}
In this paper, we implement the $3^{rd}$ order upwind-biased (Upwind3), $5^{th}$ order upwind-biased (Upwind5), and $5^{th}$ order WENO5 reconstruction methods. Upwind methods are based on solution reconstruction using the ideal weights obtained from the best polynomial approximation that matches the cell-averages in the relevant computational stencil. However, they become extremely oscillatory when the numerical solution develops discontinuities, such as shock waves, which makes the numerical solution eventually unstable. For those cases, WENO methods provide a robustness numerical solution near discontinuities while maintaining the high-order accuracy in the smooth regions of the solution by blending the ideal reconstruction weights with smoothness indicators to obtain a nonlinear reconstruction. We refer the interested readers to \cite{Shu1998} for more details.

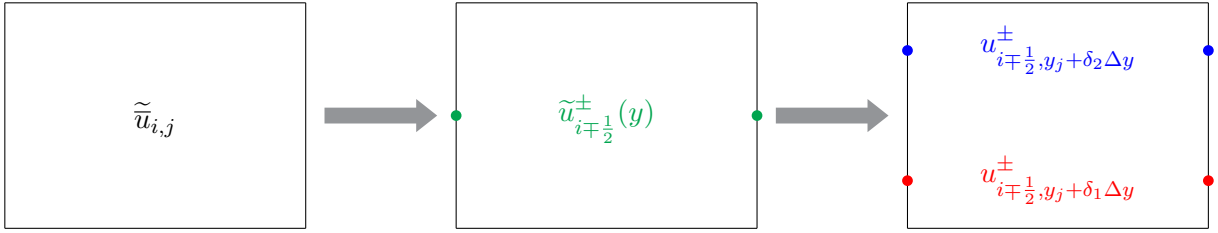
\begin{figure}[htbp]
    \begin{center}
        \begin{tikzpicture}
            \draw[-] (0,0) -- (4,0);
            \draw[-] (4,0) -- (4,3);
            \draw[-] (4,3) -- (0,3);
            \draw[-] (0,3) -- (0,0);
            \draw[-] (6,0)  -- (10,0);
            \draw[-] (10,0) -- (10,3);
            \draw[-] (10,3) -- (6,3);
            \draw[-] (6,3)  -- (6,0);
            \draw[-] (12,0) -- (16,0);
            \draw[-] (16,0) -- (16,3);
            \draw[-] (16,3) -- (12,3);
            \draw[-] (12,3) -- (12,0);
            \filldraw[very thick,Green] (6,1.5)  circle (0.05);
            \filldraw[very thick,Green] (10,1.5)  circle (0.05);
            \filldraw[very thick,red] (12,0.634)  circle (0.05);
            \filldraw[very thick,blue] (12,2.366)  circle (0.05);
            \filldraw[very thick,red] (16,0.634)  circle (0.05);
            \filldraw[very thick,blue] (16,2.366)  circle (0.05);
            \node (p0) at (2, 1.45) {$\widetilde{\overline{u}}_{i,j}$};
            \node (p1) at (8, 1.45) {\color{Green}$\widetilde{u}_{i\mp\frac{1}{2}}^\pm(y)$};
            \node (p3) at (14, 0.634) {\color{red}$u_{i\mp\frac{1}{2},y_j+\delta_1\Delta y}^\pm$};
            \node (p4) at (14, 2.366) {\color{blue}$u_{i\mp\frac{1}{2},y_j+\delta_2\Delta y}^\pm$};
            \draw[-{Triangle[width=12pt,length=8pt,Gray]}, line width=6pt,Gray](4.25,1.5) -- (5.75, 1.5);
            \draw[-{Triangle[width=12pt,length=8pt,Gray]}, line width=6pt,Gray](10.25,1.5) -- (11.75, 1.5);
        \end{tikzpicture}
        \caption{Step-by-step $3^{rd}$ order reconstruction from cell averages in the $x$-direction}\label{fig:step-by-step-recon}
    \end{center}
\end{figure}

To obtain a high-order 2-dimensional finite volume discretization, we follow the dimension-by-dimension reconstruction method in \cite{Shu1998,shu2002}. This involves two 1-dimensional reconstruction steps for each direction. In \Cref{fig:step-by-step-recon}, the reconstruction procedure for the $3^{rd}$ order reconstruction in the $x$-direction is depicted, for which we only use two Gauss quadrature points for approximating the surface integrals. Starting with the cell average $\widetilde{\overline{u}}_{i,j}$, we perform the first 1-dimensional reconstruction in the $x$-direction at $x=x_{i\pm1/2}$. This gives 1-dimensional cell averages $\widetilde{u}_{i\mp1/2}^\pm$ in the $y$-direction at $x=x_{i\pm1/2}$. This is followed by the second 1-dimensional reconstruction, now in the $y$-direction, to obtain point-wise values $u_{i\mp1/2,y_j+\delta_2\Delta y}^\pm$ at the Gauss quadrature points at $x=x_{i\pm1/2}$ and $y\in\{y_j+\delta_1\Delta y,y_j+\delta_2\Delta y\}$. 

\subsubsection{Step 1: ``De-cell'' averaging in $x$}
We now will discuss the reconstruction of $\widetilde{u}_{i+\frac{1}{2},j}^-$, which is a 1-dimensional cell average in $y$. The procedure to reconstruct $\widetilde{u}_{i-\frac{1}{2},j}^+$ is similar (see \cite{Shu1998}) therefore it will not be presented here. 
First, we set $v_{i,j}=\widetilde{\overline{u}}_{i,j}$ but the index $j$ is dropped in $v_{i,j}$ below for simplicity: 
\begin{equation}
    \widetilde{u}_{i+\frac{1}{2},j}^-=\sum_{r=0}^{k}\omega_rv^{(r)}_{i+1/2}
\end{equation}
For the $3^{rd}$ order upwind method (Upwind3), $k=1$ and 
\begin{equation}
    \begin{aligned}
        v^{(0)}_{i+1/2} &= \frac{1}{2}\left(v_{i} + v_{i+1}\right), \\
        v^{(1)}_{i+1/2} &= \frac{1}{2}\left(-v_{i-1} + 3v_{i}\right), 
    \end{aligned}
\end{equation}
while for the $5^{th}$ order Upwind5 and WENO5 schemes $k=2$ and
\begin{equation}
    \begin{aligned}
        v^{(0)}_{i+1/2} &= \frac{1}{6}\left(2v_{i} + 5v_{i+1} -  v_{i+2}\right), \\
        v^{(1)}_{i+1/2} &= \frac{1}{6}\left(-v_{i-1} + 5v_{i} + 2v_{i+1}\right), \\
        v^{(2)}_{i+1/2} &= \frac{1}{6}\left(2v_{i-2} - 7v_{i-1} + 11v_{i }\right). \\
    \end{aligned}
\end{equation}
In the upwind methods  $\omega_r$ are determined from the ideal linear weights. For example, $\omega_0=2/3$ and $\omega_1=1/3$ for Upwind3 while $\omega_0=3/10$, $\omega_1=3/5$ and $\omega_2=1/10$ for Upwind5. For the WENO5 reconstruction, we compute the nonlinear weights using 
\begin{equation}
    \omega_r = \frac{\alpha_r}{\sum_{s=0}^2\alpha_s},
\end{equation}
where $\alpha_r=d_r/\left(\beta_r+\varepsilon\right)^2$ for $r=0,1,2$, $\varepsilon=\Delta x^2$ as suggested by \cite{don2013}, $d_0=3/10$, $d_1=3/5$ and $d_2=1/10$ (same as the ideal linear weights), and the smoothness indicators are given as 
\begin{equation}
    \begin{aligned}
        \beta_0 &= \frac{13}{12}\left(v_{i} - 2v_{i+1} + v_{i+2}\right)^2 + \frac{1}{4}\left(3v_{i} - 4v_{i+1} + v_{i+2}\right)^2, \\
        \beta_1 &= \frac{13}{12}\left(v_{i-1} - 2v_{i} + v_{i+1}\right)^2 + \frac{1}{4}\left(v_{i-1}- v_{i+1}\right)^2, \\
        \beta_2 &= \frac{13}{12}\left(v_{i-2} - 2v_{i-1} + v_{i}\right)^2 + \frac{1}{4}\left(v_{i-2} - 4v_{i-1} + 3v_{i}\right)^2. \\
    \end{aligned}
\end{equation}

\subsubsection{Step 2: ``De-cell'' averaging in $y$}
Now that we have 1-dimensional cell averages in $y$, we will next complete the ``de-cell'' averaging by reconstructing point-wise values at $(x_{i\pm1/2},y_j+\delta_m\Delta y)$ for each quadrature point $m$. The procedure is similar to Step 1, but linear and nonlinear reconstruction weights are specifically defined for each quadrature point $m$.

As in Step 1, we will set $v_{i+1/2,,j}=\widetilde{u}_{i+\frac{1}{2},j}^-$ but drop the index $i+1/2$ in $v_{i+1/2,,j}$ for simplicity. Then,
\begin{equation}\label{eq:fv-recon-step2}
    {u}_{i+\frac{1}{2},y_j+\delta_m\Delta y}^-=\sum_{r=0}^{k}\omega_m^{(r)}v^{(r)}_{m},
\end{equation}
where
\begin{equation}
    v^{(r)}_{m} = \sum_{l=0}^{k} c_{rl}^mv_{j-r+l}, 
\end{equation}
and $\{c_{rl}\}^m$ denotes the components of the coefficient matrix $C^m$ for the particular quadrature point $m$. 

For the Upwind3 scheme where $k=1$, we perform the reconstruction procedure on 2 quadrature points. The coefficient matrix for each quadrature point is given as 
\begin{equation}
    \begin{aligned}
        C^{m=1}=
        \begin{pmatrix*}[r]
         808/627  & -390/1351  \\
         390/1351 &  961/1351
        \end{pmatrix*},\qquad
        C^{m=2}=
        \begin{pmatrix*}[r]
             961/1351 & 390/1351  \\
            -390/1351 & 808/627   
        \end{pmatrix*},
    \end{aligned}
\end{equation}
and the linear weights are $\omega^{(r)}_m=1/2$ for $r=0,1$ and $m=1,2$.

For the Upwind5 and WENO5 schemes where $k=2$, we consider a reconstruction on 3 quadrature points. The coefficient matrices for these points is given as 
\begin{equation}
    \begin{aligned}
        C^{m=1}&=
        \begin{pmatrix*}[r]
            4725/2927  & -2051/2438   &   249/1097  \\
             249/1097  &     14/15    &  -467/2913  \\
            -467/2913  &   366/517    &  2209/4883
        \end{pmatrix*}\\
        C^{m=2}&=
        \begin{pmatrix*}[r]
              23/24  &  1/12   & -1/24 \\    
              -1/24  & 13/12   & -1/24 \\    
              -1/24  &  1/12   & 23/24 
        \end{pmatrix*}\\
        C^{m=3}&=
        \begin{pmatrix*}[r]
            2209/4883 &   366/517   & -467/2913 \\ 
            -467/2913 &     14/15   &  249/1097 \\ 
             249/1097 & -2051/2438  & 4725/2927 
        \end{pmatrix*}.
    \end{aligned}
\end{equation}

For the $5^{th}$ order schemes with 3 quadrature points along the cell interfaces, the ideal linear coefficients of the second quadrature point, $m=2$, become negative and these are treated according to the method presented in \cite{shu2002}. Following their method, we first set the linear weights for the first and the last quadrature point as:
\begin{equation*}
    \begin{aligned}
        \gamma^{(r=0)}_{m=1}=882/6305,\quad \gamma^{(r=1)}_{m=1}=403/655,\quad \gamma^{(r=2)}_{m=1}=463/1891  \\ 
        \gamma^{(r=0)}_{m=3}=463/1891,\quad \gamma^{(r=1)}_{m=3}=403/655 , \quad \gamma^{(r=2)}_{m=3}=882/6305  \\ 
    \end{aligned}
\end{equation*}
Then, the split coefficients are used for the second quadrature point:
\begin{equation*}
    \begin{aligned}
        &\gamma^{(r=0)+}_{m=2}=9/214,\quad 
        \gamma^{(r=1)+}_{m=2}=98/107,\quad 
        \gamma^{(r=2)+}_{m=2}= 9/214  \\ 
        &\gamma^{(r=0)-}_{m=2}=\:\;9/67,\quad 
        \gamma^{(r=1)-}_{m=2}=\:\;49/67, \quad 
        \gamma^{(r=2)-}_{m=2}= \:\:\,9/67  \\ 
    \end{aligned}
\end{equation*} 

The Upwind5 scheme simply sets $\omega^{(r)}_m=\gamma^{(r)}_m$ for $m=1,3$ and $r=0,1,2$ as these coefficients are already positive, and then, it applies \Cref{eq:fv-recon-step2}. For the second quadrature point $m=2$, we perform the reconstruction using
\begin{equation}\label{eq:neg-coeff-treatment}
    {u}_{i+\frac{1}{2},y_j+\delta_2\Delta y}^-=\sigma^+u^+-\sigma^-u^-,
\end{equation}
where $\sigma^+=107/40$ and $\sigma^-=67/40$ with
\begin{equation}
    u^\pm=\sum_{r=0}^{k}\gamma^{(r)\pm}_{m=2}v^{(r)}_{m=2},
\end{equation}

The WENO5 scheme is similar to Upwind5 but it computes the nonlinear weights, instead. For $m=1,3$, it uses \Cref{eq:fv-recon-step2} with the nonlinear weights
\begin{equation}
    \omega_m^{(r)} = \frac{\alpha_m^{(r)}}{\sum_{s=0}^2\alpha_m^{(s)}},
\end{equation}
where $\alpha_m^{(r)}=\gamma_m^{(r)}/\left(\beta_r+\varepsilon\right)^2$ for $r=0,1,2$ and the smoothness indicators are given as 
\begin{equation}
    \begin{aligned}
        \beta_0 &= \frac{13}{12}\left(v_{j} - 2v_{j+1} + v_{j+2}\right)^2 + \frac{1}{4}\left(3v_{j} - 4v_{j+1} + v_{j+2}\right)^2, \\
        \beta_1 &= \frac{13}{12}\left(v_{j-1} - 2v_{j} + v_{j+1}\right)^2 + \frac{1}{4}\left(v_{j-1}- v_{j+1}\right)^2, \\
        \beta_2 &= \frac{13}{12}\left(v_{j-2} - 2v_{j-1} + v_{j}\right)^2 + \frac{1}{4}\left(v_{j-2} - 4v_{j-1} + 3v_{j}\right)^2. \\
    \end{aligned}
\end{equation}
For the second quadrature point $m=2$, WENO5 uses \Cref{eq:neg-coeff-treatment} with the same $\sigma^\pm$ but it sets $u^\pm$ as
\begin{equation}
    u^\pm=\sum_{r=0}^{k}\omega^{(r)\pm}_{m=2}v^{(r)}_{m=2},
\end{equation}
where $\alpha_{2}^{(r)\pm}=\gamma_{m=2}^{(r)\pm}/\left(\beta_r+\varepsilon\right)^2$ for $r=0,1,2$.

For the reconstruction procedure in the y-direction, the above-mentioned procedures are repeated by applying Step 1 in the y-direction to construct $\overline{u}_{i,j+1/2}$ (a 1-dimensional cell average in x) and Step 2 in the x-direction to compute $u_{x_i+\delta_m\Delta x,y_{j\mp1/2}}^\pm$ on quadrature points $m$. 
%
%
%
%
\subsection{Finite Volume Method for the Shallow Water Equations}\label{sec:fv-system-review}
In this study, we solve the cell-averaged SWEs,
\begin{equation}\label{eq:SWE-system-cell-avg}
    \frac{\partial \widetilde{\overline{\Q}}_{i,j}}{\partial t} + \frac{\widehat{\Fc}_{i+1/2,j}-\widehat{\Fc}_{i-1/2,j}}{\Delta x}+\frac{\widehat{\Gc}_{i,j+1/2}-\widehat{\Gc}_{i,j-1/2}}{\Delta y}= \widetilde{\overline{\Sc}}_{i,j},
\end{equation}
on a uniform Cartesian mesh using the finite volume method where the fluxes vectors are computed by the Gauss-Legendre quadrature rule 
\begin{equation}
    \begin{aligned}
        \widehat{\Fc}_{i\pm1/2,j}&=\sum_{m}w_m \widehat{\Fc}\left(\Q_{i\pm1/2,y_j+\delta_m\Delta y}^-,\Q_{i\pm1/2,y_j+\delta_m\Delta y}^+\right), \\
        \widehat{\Gc}_{i,j\pm1/2}&=\sum_{m}w_m \widehat{\Gc}\left(\Q_{x_i+\delta_m\Delta x,j\pm1/2}^-,\Q_{x_i+\delta_m\Delta x,j\pm1/2}^+\right),
    \end{aligned}
\end{equation}
with the Local Lax-Friedrichs flux
\begin{equation}
    \begin{aligned}
        \widehat{\Fc}\left(\Q^-,\Q^+\right) &=  \frac{1}{2}\left(\Fc(\Q^-)+\Fc(\Q^+)\right)-\frac{\lambda_F}{2}\left(\Q^+-\Q^-\right),\\
        \widehat{\Gc}\left(\Q^-,\Q^+\right) &=\frac{1}{2}\left(\Gc(\Q^-)+\Gc(\Q^+)\right)-\frac{\lambda_G}{2}\left(\Q^+-\Q^-\right),
    \end{aligned}
\end{equation}
where $\lambda_F$ and $\lambda_G$ are the maximum eigenvalues of the flux Jacobians $\left|{\partial \Fc}/{\partial \Q}\right|$ and $\left|{\partial \Gc}/{\partial \Q}\right|$, respectively.

The high-order reconstruction is performed by applying the procedures discussed in \Cref{sec:fv-recon} to $\Q$ in a component-by-component fashion.
%
%
%
%
\section{Tensor Train Decomposition}\label{sec:TT-decomposition}

In this section we briefly introduce the tensor notation and the tensor train manipulation techniques we apply in this work.

\subsection{Tensor Train}
Tensor train, or TT-format of a tensor, represents the tensor as a product of cores, which are either matrices or three-dimensional tensors~\cite{oseledets2010tt}. Generally, a TT representation $\mathcal{X}^{TT}$ of a d-dimensional tensor $\mathcal{X}$ is defined as
\begin{equation}
    \label{eq:TT_def}
    \ten{X}_{TT}(i_1,\ldots,i_d) = \sum_{\alpha_1,\ldots,\alpha_{d-1}}^{r_1,\ldots,r_{d-1}}\ten{G}_1(1,i_1,\alpha_1)\ldots\ten{G}_d(\alpha_{d-1},i_d,1) + \varepsilon(i_1,\ldots,i_d),
\end{equation}
where the error, $\varepsilon$, is a tensor with the same dimensions as $\ten{X}$. The integers $[r_1,\ldots, r_{d-1}]$ are called TT-ranks, and $\ten{G}k$ are called TT-cores.
The last term, $\varepsilon$, is a tensor with the same
dimensions of $\ten{X}$ representing the approximation error.
Equivalently, we can also denote the TT-format by the multiple
matrix product,

\begin{equation}
  \ten{X}_{TT}(i_1,i_2,\dots,i_d) =
  \mat{G}_1(i_1)\mat{G}_2(i_2)\dots\mat{G}_d(i_d) + \varepsilon,
  \label{eqn:TT-vector}
\end{equation}
where each term
$\big(\mat{G}_{k}(i_k)\big)_{\alpha_{k-1},\alpha_{k}}$,
$i_k=1,2,\ldots,n_k$, $k=1,2,\ldots,d$, is a matrix of size
$r_{k-1}\times r_k$.
Therefore, the $\ten{G}_k(:,i_k,:)$ are a set of matrix
slices $\mat{G}_k(i_k)$ that are labeled with the single index $i_k$.
The entries of the integer array
$\mathbf{r}=\big[r_1,\dots,r_{d-1}\big]$ are the TT-ranks, and
quantify the compression effectiveness.
Since each TT-core only depends on a single mode index of the full
tensor $\ten{X}$, e.g., $i_k$, the TT-format effectively embodies a
discrete separation of variables.
When the TT-ranks are relatively small with respect to the problem
size, a TT-based approach is referred to as a \emph{low-rank
approximation} \cite{bachmayr2023low}.

Assuming that $n_k=\mathcal{O}(n)$ and $r_k=\mathcal{O}(r)$ for some
non-negative integers $n$ and $r$, and for all $k=1,2,\ldots,d$, the
total number of elements that TT-format stores is proportional to
$\ten{O}(2nr+(d-2)nr^2)$, which is linear with the number of
dimensions $d$.

In the case that $\ten{X}$ is a matrix, its TT-format is simplified to the following,
\begin{equation}
    \label{eq:TT_def_2D}
    \ten{X}_{TT}(i_1,i_2) =\mat{G}_1(i_1)\mat{G}_2(i_2) + \varepsilon(i_1,i_2)
\end{equation}
Given that the shallow water problems we investigate in this work have two spatial dimensions, this decomposition will be the one we use to represent the solutions at each time step.
\subsection{TT Rounding}
\label{rounding}
When a matrix is represented in TT-format as $\ten{X}_{TT}$ with a TT-rank of $r_1$, the TT rounding procedure is used to obtain a more compact TT representation, $\ten{Y}_{TT}$, with TT-ranks $r'_1 \leq r_1$, while maintaining a specified relative accuracy $\varepsilon_{TT}$. This procedure is also called truncation or recompression~\cite{oseledets2011tensor}.
For matrices, the TT rounding algorithm consists of two steps. First, the second core $\mat{G}_2$ is orthogonalized using $RQ$ decomposition. Thereafter, singular value decomposition (SVD) truncate at tolerance $\varepsilon_{TT}$ is applied to the product $\mat{G}_1\mat{R}$ to arrive at a new decomposition $\ten{Y}_{TT}$. In this study, we denote TT rounding by $\ten{Y}_{TT}=\text{round}(\ten{X}_{TT},\varepsilon_{TT})$.

\subsection{TT Cross Interpolation}

TT cross interpolation is a technique used to construct a TT representation of a tensor without needing to form the entire tensor explicitly. This method is particularly valuable when dealing with very large tensors or in situations where calculations are impractical due to limitations in TT arithmetic.
Stemming from the skeleton (or CUR) decomposition~\cite{mahoney2009cur}, in combination with the Maximum Volume Principle~\cite{goreinov2010find}, heuristic cross interpolation algorithms for tensor train, such as Alternating Minimal Energy (AMEn) \cite{dolgov2014alternating}, or Density Matrix Renormalization Group (DMRG) \cite{savostyanov2011fast} have been developed. In this work, we use an implementation of AMEn algorithm, \texttt{amen\_cross}, which is available in MATLAB TT-Toolbox~\cite{tt-toolbox}.


\section{Tensorization of the Finite Volume Scheme}\label{sec:tensorization-of-the-method}
In this section, we will discuss the tensorization of the finite volume method for the shallow water equations. 
\subsection{Tensor Train Finite Volume (TT-FV) method for hyperbolic conservation laws}
Following \cite{danis2024tensortrain}, we start with the full-tensor form of the SWE equations with ``\emph{loop indices}''
\begin{equation}\label{eq:SWE-system-for-loop}
    \frac{\partial \widetilde{\overline{\Q}}_{i,j}}{\partial t} + \frac{\widehat{\Fc}_{i+1/2,j}-\widehat{\Fc}_{i-1/2,j}}{\Delta x}+\frac{\widehat{\Gc}_{i,j+1/2}-\widehat{\Gc}_{i,j-1/2}}{\Delta y}= \widetilde{\overline{\Sc}}_{i,j}.
\end{equation}

On a structured Cartesian mesh, we can introduce the shift operators in $x$ and $y$, and rewrite the flux terms as
\begin{equation}\label{eq:shifted-fluxes}
    \begin{aligned}
        \Fchat_{i+\frac{1}{2},j} &= T^x_{i,j}\Fchat_{i-\frac{1}{2},j}, \\
        \Gchat_{i,j+\frac{1}{2}} &= T^y_{i,j}\Gchat_{i,j-\frac{1}{2}}.
    \end{aligned}
\end{equation}
Substituting \Cref{eq:shifted-fluxes} into \Cref{eq:SWE-system-for-loop} and dropping the \emph{loop indices}, we obtain the vectorized form of the SWEs, which we will refer to as the ``\emph{full-tensor}'' form:
\begin{equation}\label{eq:SWE-system-no-loop}
    \frac{\partial \widetilde{\overline{\Q}}}{\partial t} + \frac{1}{\Delta x}\left(T^x-1\right)\widehat{\Fc}+\frac{1}{\Delta y}\left(T^y-1\right)\widehat{\Gc}= \widetilde{\overline{\Sc}}.
\end{equation}
Note that terms in \Cref{eq:SWE-system-no-loop} correspond to 2-dimensional pre-stored arrays. Therefore, we simply replace the full-tensor terms with their TT counterparts:
\begin{equation}\label{eq:SWE-system-TT}
    \frac{\partial \widetilde{\overline{\Q}}_{TT}}{\partial t} + \frac{1}{\Delta x}\left(T^x-1\right)\widehat{\Fc}_{TT}+\frac{1}{\Delta y}\left(T^y-1\right)\widehat{\Gc}_{TT}= \widetilde{\overline{\Sc}}_{TT}.
\end{equation}
\subsection{Computing the fluxes in the tensor-train format}
In this study, we implement the TT format of the Local Lax-Friedrichs flux similar to that suggested by \cite{mustafa2023sod}. For example, the fluxes in the $x$-direction will be computed in the TT format as
\begin{equation}\label{eq:LLF-TT}
    \widehat{\Fc}_{TT}\left(\Q^-_{TT},\Q^+_{TT}\right) = \frac{1}{2}\left(\Fc(\Q^-_{TT})+\Fc(\Q^+_{TT})\right)-\frac{\lambda_{F,TT}}{2}\left(\Q^+_{TT}-\Q^-_{TT}\right)
\end{equation}
However, the linear and nonlinear SWEs differ in the computation of each individual term in \Cref{eq:LLF-TT}. For the linear SWEs the physical flux terms, $\Fc(\Q^\pm_{TT})$ can be directly computed, without relying on special considerations such as TT cross interpolation. In addition, the eigenvalue $\lambda_{F}=\sqrt{gH}$ is a constant for the linear equations. In contrast, for the nonlinear SWE equations $1/h_{TT}$ must be computed for the physical fluxes $\Fc$ and $\Gc$ and $\lambda_{F,TT}=|u|+\sqrt{gh_{TT}}$ must be computed as the eigenvalue of the flux Jacobian. None of these quantities can be computed directly, therefore we employ the AMEn method to compute eigenvalues $\lambda_{F,TT}=|u|+\sqrt{gh_{TT}}$ and the Taylor series approximation given in \Cref{alg:taylor-inverse}, to compute the reciprocal of the tensor train $h_{TT}$: 
\begin{algorithm}[htbp]
    \caption{Taylor Series Approximation to Tensor Train Reciprocal}\label{alg:taylor-inverse}
    \KwData{$x_{TT}$, $\varepsilon_{TT}>0$}
    \KwResult{$y_{TT}=1/x_{TT}$}
    $y_{TT} \gets 1$\; 
    $\Delta y_{TT} \gets 1$\; 
    $Err \gets 1$\; 
    $N \gets \text{numel}(x_{TT})$\;
    $x_{avg} \gets \text{sum}(x_{TT})/N$\; 
    $\tilde{x}_{TT} \gets \text{round}(1-x_{TT}/x_{avg},\varepsilon_{TT})$\; 
    \While{$Err > \varepsilon_{TT}$}{
        $\Delta y_{TT} \gets \text{round}(\Delta y_{TT}*\tilde{x}_{TT},\varepsilon_{TT})$\; 
        $y_{TT} \gets \text{round}(y_{TT} + \Delta y_{TT},\varepsilon_{TT})$\;
        $Err \gets \|\Delta y_{TT}\|_F/\sqrt{N}$\;
    }
    $y_{TT} \gets y_{TT}/x_{avg}$;
\end{algorithm}

After obtaining $1/h_{TT}$, the components of the flux vectors in \Cref{eq:nonlinear-SWE} can be directly computed. Note that we could have also employed a TT-cross interpolation method to compute $1/h_{TT}$ instead of \Cref{alg:taylor-inverse}. In the numerical examples considered in this study, however, we found that Taylor series approximation to $1/h_{TT}$ is as fast as the AMEn method. This is possibly because the shallow water layer thickness $h$ can be decomposed as a superposition of a large mean value and small amplitude oscillations, i.e. $h=H+\eta(x,y,t)$ where $|\eta|\ll H$, which leads to a very fast and robust convergence of the Taylor series approximation for $1/h_{TT}$.

Note also that this method is different than the LF-cross method developed in \cite{danis2024tensortrain}, where the complete flux vector of the compressible Euler equations is computed with a single cross interpolation using the AMEn method. In the context of the finite difference method for solving compressible flow equations, the LF-cross method was reported to be faster than a similar method presented here that approximates $1/\rho_{TT}$ (reciprocal of density) with the TT cross interpolation to compute the flux vector components of compressible Euler equations directly. This was thought to be due to slow TT-rounding while each flux vector component was estimated. However, for the finite volume method for solving the shallow water equations, we found that the LF-cross approach is considerably slower than the present approach. This might be due to two reasons. First, the compressible Euler equations have more conserved variables than the SWEs, and therefore, more floating point operations are needed to compute fluxes in compressible flow, meaning that the rounding routine will be called more often and each call will more expensive. Second, computing $1/h_{TT}$ for the SWEs is likely a much simpler task than computing the $1/\rho_{TT}$ for the compressible flow. As mentioned above, in the case of the SWEs, $h$ is simply a superposition of a mean depth and small-amplitude waves, but in the compressible flow, the density field $\rho$ can have large variations, and even, discontinuities such as strong shock waves. Therefore, the nonlinear flux calculations should be designed or selected according to the particular application, especially when the computation of the reciprocal of a tensor train is required.
\subsection{High-order reconstruction in the tensor-train format}\label{sec:fv-recon-TT}
In this paper, we consider two types of high-order variable reconstruction methods -- linear (Upwind3 and Upwind5) and nonlinear (WENO5) reconstructions as previously discussed in \Cref{sec:fv-recon}. For linear reconstruction schemes, we directly manipulate the TT cores. However, for nonlinear WENO reconstruction, this is not possible due to TT reciprocals. Therefore, we employ the TT cross interpolation for the WENO reconstruction.
\subsubsection{Linear Reconstruction}\label{sec:fv-recon-TT-linear}
The linear reconstruction in the TT format in a given direction is applied only to the TT core corresponding to that direction, which significantly reduces the computational costs. To exploit this, we modify the algorithm presented in \Cref{sec:fv-recon}. Note that, in the full-tensor format, Step 1 is followed by Step 2 for each time a reconstruction is needed along a cell interface. This means that, in a 2-dimensional setting, Step 2 is applied twice to close the cell boundaries. However, this is not needed in the TT format. 

To illustrate the idea, we adopt the notation used in \cite{oseledet2011tt} and denote the elementwise values of a cell-averaged tensor $\widetilde{\overline{u}}$ as
\begin{equation}
    \widetilde{\overline{u}}(i,j) = \overline{u}_1(i)\widetilde{u}_2(j),
\end{equation}
Note that the first core $\overline{u}_1$ is written with a bar to denote the cell-averaging operator in $x$ and the second core $\widetilde{u}_2$ is written with a tilde to denote the cell-averaging operator in $y$. This implies that Step 2 in \Cref{sec:fv-recon} can be applied to each core independently and even simultaneously.  

\emph{TT-Reconstruction Step 1} starts with reconstructing the cores at the quadrature points $u_1(x_i+\delta_m\Delta x)$ and $u_2(y_j+\delta_m\Delta y)$. This step is almost identical to Step 2 in \Cref{sec:fv-recon} except for that fact that the reconstruction is only applied to TT cores rather than to the full-tensor.

\emph{TT-Reconstruction Step 2} is similar to Step 1 in \Cref{sec:fv-recon}, and again, we only apply reconstruction to the TT-cores. For example, a reconstruction in the $x$-direction first compute $u_1(i\pm1/2)^\mp$, and then, the result is combined with $u_2(y_j+\delta_m\Delta y)$ prepared in TT-Reconstruction Step 1:
\begin{equation}
    u(i\pm1/2,y_j+\delta_m\Delta y)^\mp = u_1(i\pm1/2)^\mp u_2(y_j+\delta_m\Delta y).
\end{equation}
Similarly, for the reconstruction in the $y$-direction, TT-Reconstruction Step 2 computes $u_2(j\pm1/2)^\mp$ and combines this with $u_1(x_i+\delta_m\Delta x)$ prepared in TT-Reconstruction Step 1:
\begin{equation}
    u(x_i+\delta_m\Delta x,j\pm1/2)^\mp = u_1(x_i+\delta_m\Delta x) u_2(j\pm1/2)^\mp.
\end{equation}

\subsubsection{Nonlinear Reconstruction}\label{sec:fv-recon-TT-nonlinear}
The WENO scheme performs the nonlinear reconstruction using the TT cross interpolation by applying Step 1 and Step 2 of \Cref{sec:fv-recon} in the same order. In fact, TT-WENO Step 1 is similar to the finite difference WENO-cross method proposed in \cite{danis2024tensortrain}, i.e. Step 1 is applied in an equation-by-equation fashion using the TT cross interpolation as detailed in \Cref{alg:weno-cross-step-1}. 
\begin{algorithm}[htbp]
    \caption{TT-WENO Step 1 for the ``de-cell'' averaging in $x$}\label{alg:weno-cross-step-1}
    \KwData{A component of the cell-averaged conserved variable vector $\widetilde{\overline{v}}_{TT}$, function \emph{funWENO} and the convergence criterion $\varepsilon_{TT}$.}
    \KwResult{Cell-averaged $\widetilde{v}^\pm_{TT}$ in $y$ at the interface locations.}
    Collect the 5-cell stencil of $\widetilde{\overline{v}}_{TT}$ in the $x$-direction into the array $\mathbf{S} = \{T_x^{-2}\widetilde{\overline{v}}_{TT},T_x^{-1}\widetilde{\overline{v}}_{TT},\widetilde{\overline{v}}_{TT},T_x^{1}\widetilde{\overline{v}}_{TT},T_x^{2}\widetilde{\overline{v}}_{TT}\}$. \\
    Set the initial guess $v_0=\widetilde{\overline{v}}_{TT}$. \\
    Perform cross interpolation for $+$ side: $\widetilde{v}^+_{TT}=$AMEn(\emph{funWENO}, $\mathbf{S}$, $v_0$, $\varepsilon_{TT}$,$+1$) \\
    Perform cross interpolation for $-$ side: $\widetilde{v}^-_{TT}=$AMEn(\emph{funWENO}, $\mathbf{S}$, $v_0$, $\varepsilon_{TT}$,$-1$) 
\end{algorithm}

Similarly, TT-WENO Step 2 applies Step 2 of \Cref{sec:fv-recon} for each conserved variable using the TT cross interpolation, see \Cref{alg:weno-cross-step-2}. However, it performs the WENO reconstruction at all quadrature points with a single cross interpolation. Note that the \emph{AMEn-cross} routine in the MATLAB Toolbox returns a TT of size $1\times N_x \times N_y \times N_q$, where $N_q$ is the number of quadrature points. Therefore, as a last step, it is reshaped and permuted to obtain a new TT of size $1\times N_x \times N_y N_q\times 1$ for the reconstruction in the $x$-direction and $1\times N_x N_q\times N_y \times 1$ for the reconstruction in the $y$-direction.
\begin{algorithm}[htbp]
    \caption{TT-WENO Step 2 for for the ``de-cell'' averaging in $y$}\label{alg:weno-cross-step-2}
    \KwData{Cell-averaged $\widetilde{v}^\pm_{TT}$ in $y$ obtained from TT-WENO Step 1, function \emph{funWENOquad} and the convergence criterion $\varepsilon_{TT}$.}
    \KwResult{$v^\pm_{TT}$ at the interface locations.}
    Collect the 5-cell stencil of $v_{TT}$ in the $x$-direction into the arrays $\mathbf{S}^\pm = \{T_y^{-2}\widetilde{v}^\pm_{TT},T_y^{-1}\widetilde{v}^\pm_{TT},v^\pm_{TT},T_y^{1}\widetilde{v}^\pm_{TT},T_y^{2}\widetilde{v}^\pm_{TT}\}$. \\
    Set the initial guesses as $v^\pm_0=\widetilde{v}^\pm_{TT}$. \\
    Perform cross interpolation for $+$ side: $v^+_{TT}=$AMEn(\emph{funWENOquad}, $\mathbf{S}^+$, $v^+_0$, $\varepsilon_{TT}$,$+1$) \\
    Perform cross interpolation for $-$ side: $v^-_{TT}=$AMEn(\emph{funWENOquad}, $\mathbf{S}^-$, $v^-_0$, $\varepsilon_{TT}$,$-1$) \\
    Reshape and permute $v^\pm_{TT}$ to return two TTs of size $1\times N_x \times N_yN_q \times 1$. 
\end{algorithm}

\subsection{Choosing a suitable \texorpdfstring{$\varepsilon_{TT}$}{Lg}}
The choice of $\varepsilon_{TT}$ determines the overall accuracy and the speed of the TT solver, as demonstrated in \cite{danis2024tensortrain}. In this study, we modify the $\varepsilon_{TT}$ formula slightly to accommodate both the $3^{rd}$ and $5^{th}$ order schemes as 
\begin{equation}\label{eq:eps-tt-formula}
    \varepsilon_{TT}=C_\varepsilon\frac{ V^{1/2} \Delta x^{p-1/2}}{\max_{q\in\Q_{TT}}{\|q\|_{F}}},
\end{equation}
where $p$ denotes the order of accuracy of the TT-FV scheme, $V$ is the volume of the computational domain, and $C_\varepsilon$ is a problem dependent variable (see \cite{danis2024tensortrain}). 

However, if one of the conserved variables is exactly zero, choosing $\varepsilon_{TT}$ according to formula may still result in rank growth. This is primarily because $\varepsilon_{TT}$ becomes comparable to the numerical noise generated for that conserved variable, which is known to result in rank growth. A robust workaround is obtained by replacing the max operator in the denominator of \Cref{eq:eps-tt-formula} with the norm of the relevant conserved variable itself. For example, consider the conserved variable $u_{TT}\in\Q_{TT}$ of the linear SWEs \Cref{eq:linear-SWE}, for which we calculate
\begin{equation}\label{eq:modified-eps-tt-formula}
    \varepsilon^{u}_{TT}=\min{\left(10^{-3},C_\varepsilon\frac{ V^{1/2} \Delta x^{p-1/2}}{{\|u\|_{F}}}\right)}
\end{equation}
at the beginning of a time step. Note that the new dynamic formula is clipped at $10^{-3}$ to avoid accuracy loss if $\|u\|_{F}\rightarrow0$. 

Equipped with \Cref{eq:modified-eps-tt-formula}, the $3^{rd}$ order strong stability-preserving Runge-Kutta method \cite{shu1988ssp,Gottlieb2001} for a given semi-discrete form of $u_{TT}$,
\begin{equation}
    \frac{du_{TT}}{dt}=L(u_{TT}),
\end{equation}
is then defined in TT format as
\begin{equation}\label{eq:tt-ssprk3}
    \begin{aligned}
        u^{(1)}_{TT}  &= \mathcal{F}(u^{n}_{TT}),\\
        u^{(2)}_{TT}  &= \text{round}\left(\frac{3}{4}u^{n}_{TT}+\frac{1}{4}\mathcal{F}(u^{(1)}_{TT}),\,\varepsilon^{u}_{TT}\right),\\
        u^{n+1}_{TT}  &= \text{round}\left(\frac{1}{3}u^{n}_{TT}+\frac{2}{3}\mathcal{F}(u^{(2)}_{TT}),\,\varepsilon^{u}_{TT}\right),\\
    \end{aligned}
\end{equation}
where $\mathcal{F}(u)$ is the forward Euler step
\begin{equation}\label{eq:tt-forward-euler}
    \mathcal{F}(u) = \text{round}\left(u+\Delta t L(u),\,\varepsilon^{u}_{TT}\right)
\end{equation}
Note that the same procedure is applied for the time integration of other conserved variables. For all other rounding operations, we use the standard formula given in \Cref{eq:eps-tt-formula}.

%
%
%
%
\section{Numerical Results}\label{sec:results}
In this section, we investigate the performance of the TT finite volume solver in the numerical examples taken from the validation suite introduced in \cite{bishnu2024}. In all examples, the nondimensional forms of \Cref{eq:SWE-system,eq:linear-SWE,eq:nonlinear-SWE} are solved (see in Appendix) on a square computational domain $\Omega=[L,L]$ using the $3^{rd}$ order strong stability-preserving Runge-Kutta \cite{shu1988ssp,Gottlieb2001}. For the $5^{th}$ order schemes, we set the time step $\Delta t$ proportional to $\Delta x^{5/3}$ to maintain accuracy and the proportionality constant is chosen such that the numerical solution remains stable for all time steps. Note that all time step sizes will be reported in terms of the nondimensional variables. In all TT simulations, \Cref{eq:eps-tt-formula,eq:modified-eps-tt-formula} are calculated using only the nondimensional variables and the $L_2-$error is defined relative to the $L_2-$error of the same reconstruction method at the coarsest grid level. Unless otherwise stated, we set $g=10\;m/s^2$, $f=10^{-4}\;1/s$, $H=1000\;m$, $c=\sqrt{gH}=100\; m/s$ and $R=c/f=10^6\;m$ as in \cite{bishnu2024}. All surface integrals are computed using Gauss-Legendre quadrature rule in the TT format, as suggested in \cite{alexandrov2023challenging}. See \Cref{tab:test-case-summary} for a summary of the numerical examples along with nondimensional time step to grid size ratios and $C_\varepsilon$ in \Cref{eq:modified-eps-tt-formula}. 
\begin{table}[h]
\begin{center}
\begin{tabular}{|l|rc|rc|rc|}
\hline
\multicolumn{1}{|c|}{\multirow{2}{*}{\bf Test Case}} & \multicolumn{2}{c|}{\bf Upwind3}                                                   & \multicolumn{2}{c|}{\bf Upwind5}                                                          & \multicolumn{2}{c|}{\bf WENO5}                                                           \\ \cline{2-7} 
\multicolumn{1}{|c|}{}                           & \multicolumn{1}{c}{$\Delta t/\Delta x$} & \multicolumn{1}{c|}{$C_\varepsilon$} & \multicolumn{1}{c}{$\Delta t/\Delta  x^{5/3}$} & \multicolumn{1}{c|}{$C_\varepsilon$} & \multicolumn{1}{c}{$\Delta t/\Delta x^{5/3}$} & \multicolumn{1}{c|}{$C_\varepsilon$} \\ \hline
Coastal Kelvin Wave   &  $5\times 10^{-5}$   	& $1$	    & $2.5\times 10^{-4}$  & $1$	& $2.5\times 10^{-4}$  & $1000$                                \\
Inertia-Gravity Wave  &  $10^{-4}$	            & $1$	    & $10^{-3}$            & $1$	& $10^{-3}$            & $500$                                 \\
Barotropic Tide       &  $2.5\times 10^{-4}$	    & $1$	    & $5\times 10^{-4}$    & $1$	& $5\times 10^{-4}$    & $1$                                   \\
Manufactured Solution &  $5\times 10^{-5}$	    & $10^{-4}$	& $5\times 10^{-3}$    & $1$	& $5\times 10^{-3}$    & $1$                                   \\ \hline
\end{tabular}
\caption{Summary of test cases}\label{tab:test-case-summary}
\end{center}
\end{table}

\subsection{Coastal Kelvin Wave}
\begin{figure}[htbp]
    \begin{center}
        \begin{subfigure}{0.49\textwidth}
            \centering
            \includegraphics[width=1\textwidth,trim={1cm 6cm 1cm 6cm},clip]{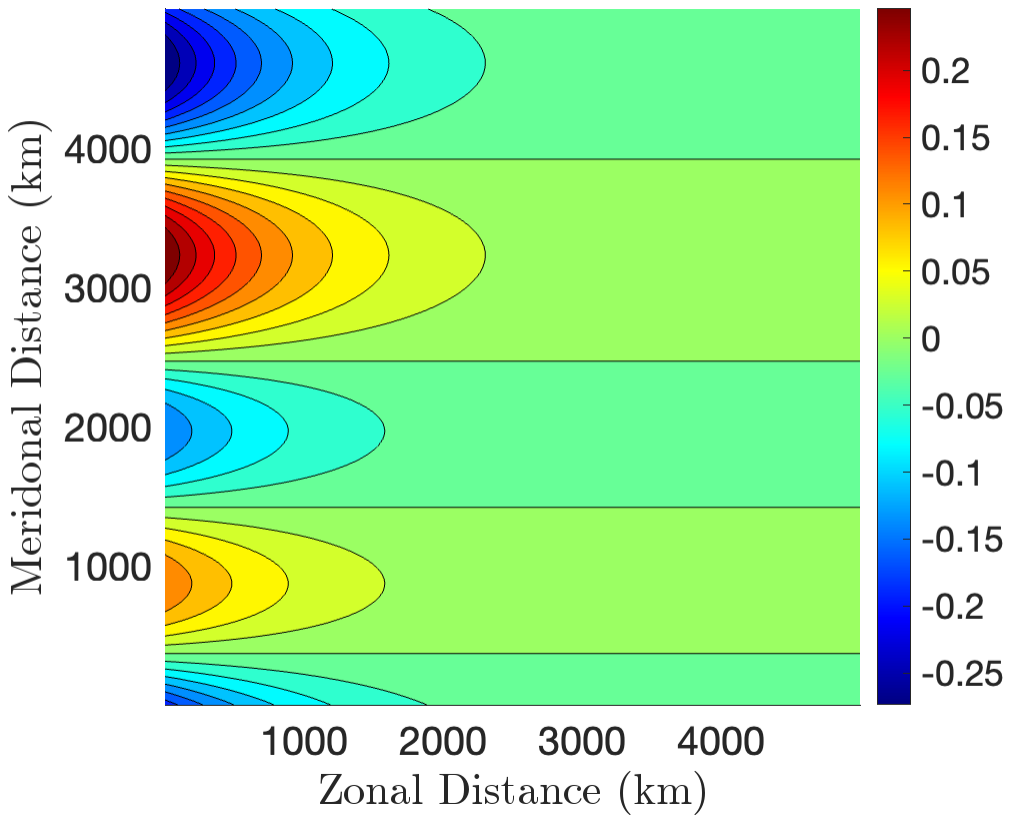}
            \caption{Contour plot of $\eta$}
        \end{subfigure}
        \begin{subfigure}{0.49\textwidth}
            \centering
            \includegraphics[width=0.95\textwidth,trim={1.2cm 6cm 1.2cm 6cm},clip]{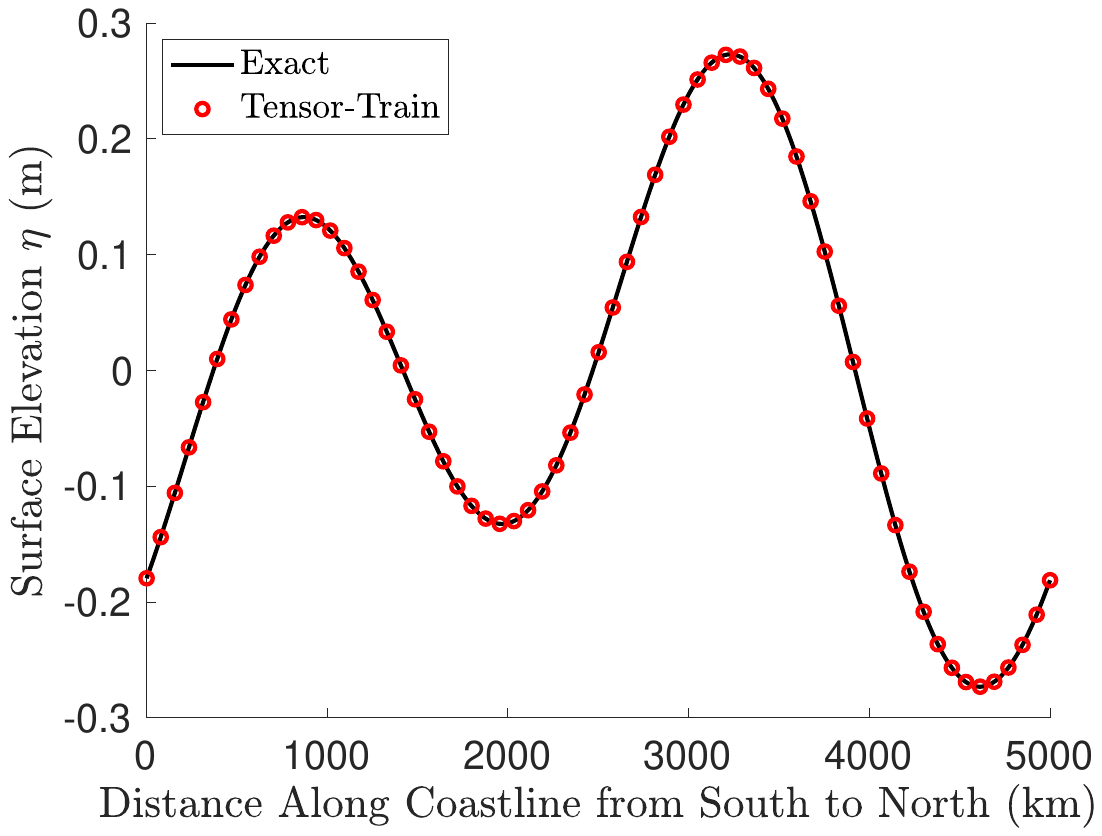}
            \caption{$\eta$ profile}
        \end{subfigure}
        \caption{Exact and TT solutions of surface elevation for Coastal Kelvin Wave at $T=3$ hours. The TT solution is obtained with Upwind5 on $1280\times1280$ grid.}\label{fig:coastal-Kelvin-profiles}
    \end{center}
\end{figure}

In this example, we solve the linear SWEs, \Cref{eq:SWE-system,eq:linear-SWE} up to the final time $T=3$ hours  and set $L=5\times10^6\;m$. Coastal Kelvin Waves occur when the Coriolis force is balanced against a topographic boundary, and are found in real-world observations\cite{JOHNSON199029}. In this idealized configuration, the coastline is in the non-periodic $x$ direction and nondispersive waves travel along the boundary in the periodic $y$ direction. The exact solution is plotted in \Cref{fig:coastal-Kelvin-profiles} and the initial conditions are composed of two wave modes:
\begin{equation}
    \begin{aligned}
        \eta&=-H\left\{\hat{\eta}^{(1)}\sin{\left(k_y^{(1)}(y+ct)\right)}+\hat{\eta}^{(2)}\sin{\left(k_y^{(1)}(y+ct)\right)}\right\}\exp{\left(-x/R\right)}, \\
        u&=0, \\
        v&=\sqrt{gH}\left\{\hat{\eta}^{(1)}\sin{\left(k_y^{(1)}(y+ct)\right)}+\hat{\eta}^{(2)}\sin{\left(k_y^{(1)}(y+ct)\right)}\right\}\exp{\left(-x/R\right)},
    \end{aligned}
\end{equation}
where $\hat{\eta}^{(1)}=\hat{\eta}^{(2)}/2=10^{-4}$ m, and $k^{(1)}_y=k^{(2)}_y/2=2\pi/L$.

We use the exact solution to set the boundary condition in the $x-$direction, but in the $y-$direction, we consider a periodic boundary. Furthermore, for both upwind methods we use $C_\varepsilon=1$ in \Cref{eq:eps-tt-formula,eq:modified-eps-tt-formula} to calculate $\varepsilon_{TT}$, but for WENO5 we uses $C_\varepsilon=1000$. To guarantee accuracy and stability of the numerical solution, we also consider a time step of $\Delta t/\Delta x=5\times10^{-5}$ for the Upwind3 method and $\Delta t/\Delta x^{5/3}=2.5\times10^{-4}$ for the $5^{th}$ order methods.

\begin{figure}[htbp]
    \begin{center}
        \begin{subfigure}{0.49\textwidth}
            \centering
            \includegraphics[width=1\textwidth,trim={1cm 6cm 1cm 6cm},clip]{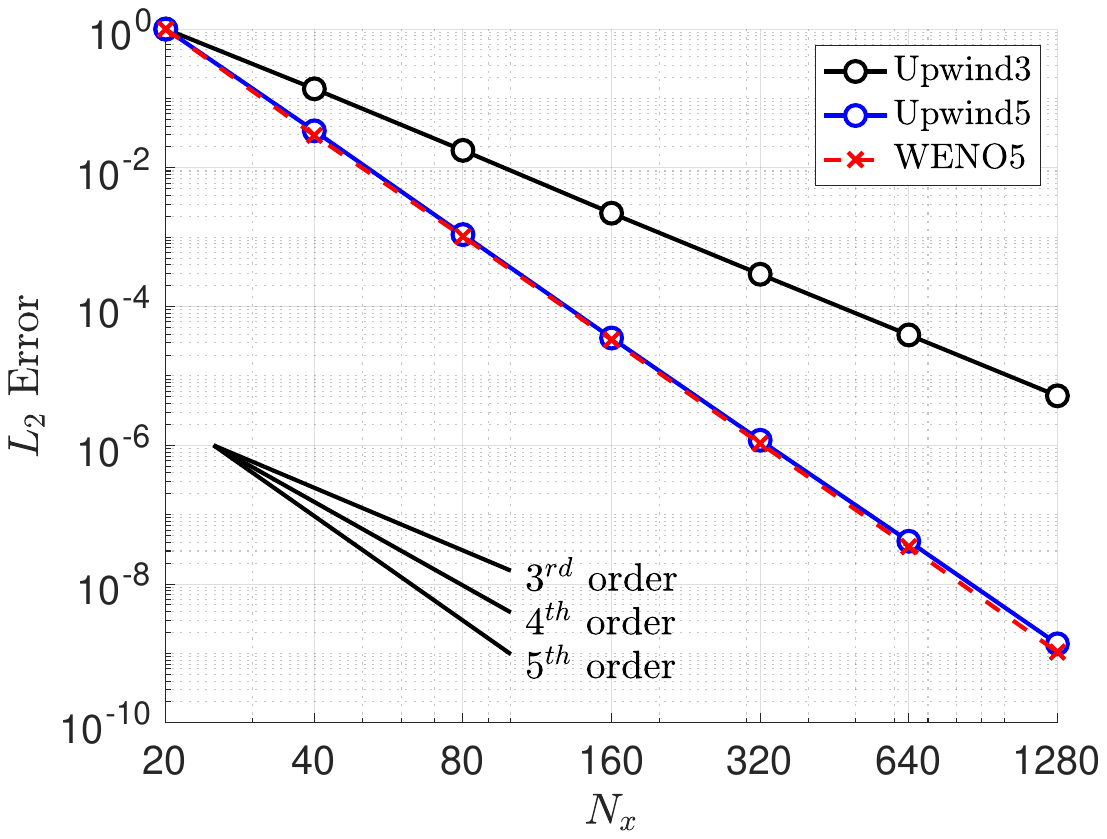}
            \caption{$L_2$ error for $\eta$}
        \end{subfigure}
        \begin{subfigure}{0.49\textwidth}
            \centering
            \includegraphics[width=0.95\textwidth,trim={1.2cm 6cm 1.2cm 6cm},clip]{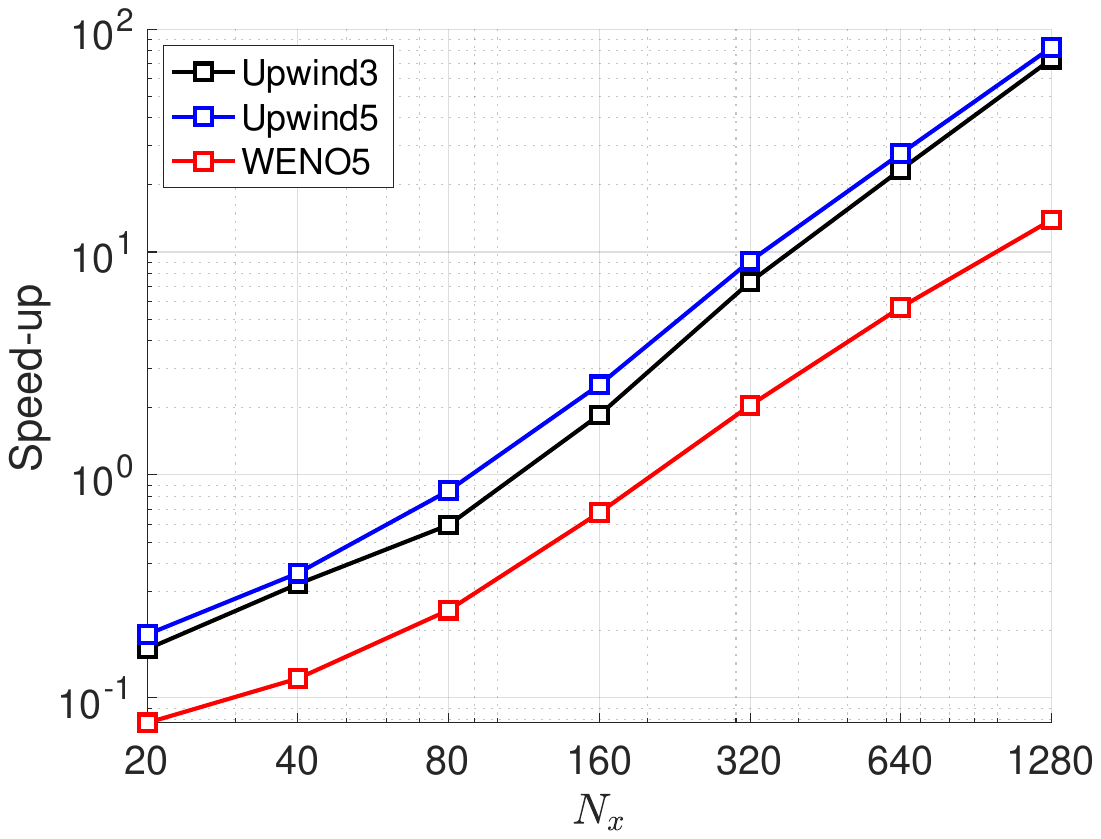}
            \caption{Speed-up}
        \end{subfigure}
        \caption{Speed-up and $L_2$ error of surface elevation for Coastal Kelvin Wave at $T=3$ hours}\label{fig:coastal-Kelvin-performance}
    \end{center}
\end{figure}

\Cref{fig:coastal-Kelvin-performance} shows the performance of Upwind3, Upwind5 and WENO5 methods. On the left, all TT reconstruction methods are observed to maintain their formal $L_2-$convergence order for $\eta$. On the right, the various methods are compared in terms of the acceleration with respect to their corresponding standard full-tensor versions: At the finest grid level, TT-Upwind5 achieves 83x acceleration while Upwind3 has a speed-up of 73x. Note that Upwind3 is still less expensive than Upwind5. Therefore, the higher speed-up value of Upwind5 simply means that it becomes more efficient in accelerating its full-tensor version for this linear problem. On the other hand, the TT-WENO5 method only achieves the significantly lower speed-up of 14x. This clearly shows the cost of the TT cross interpolation used in the WENO5 scheme, even when the problem is linear and smooth.
\subsection{Inertia-Gravity Wave}
\begin{figure}[htbp]
    \begin{center}
        \begin{subfigure}{0.49\textwidth}
            \centering
            \includegraphics[width=1\textwidth,trim={1cm 6cm 1cm 6cm},clip]{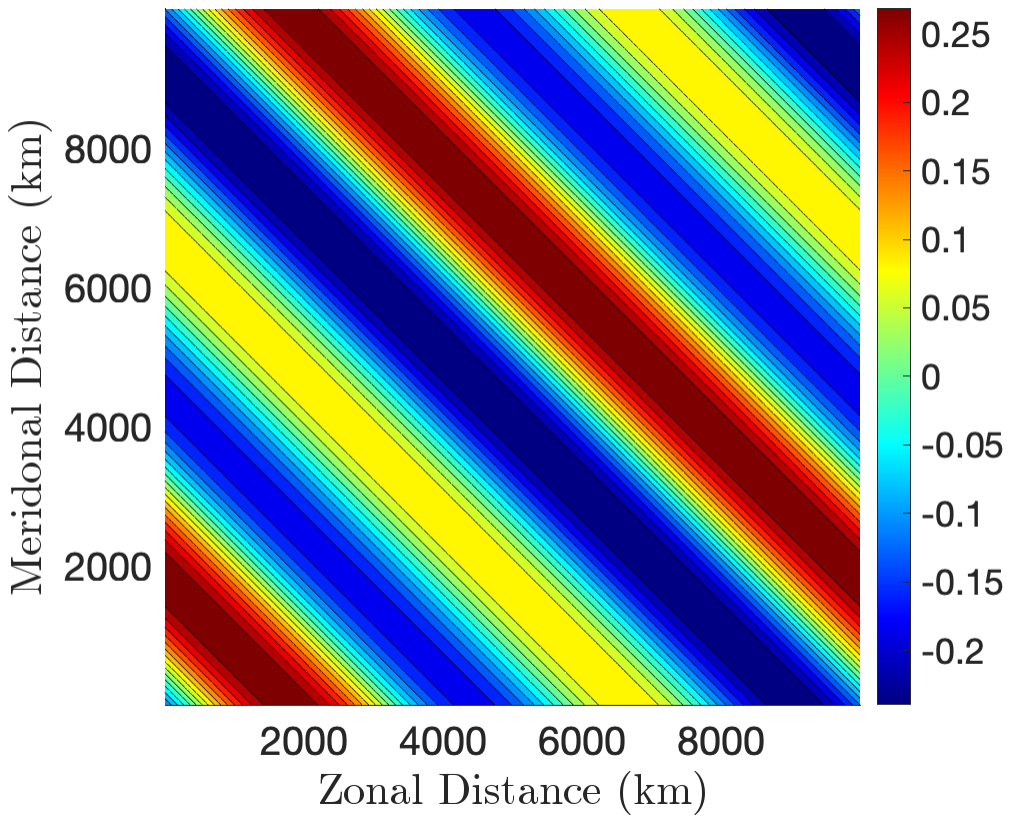}
            \caption{Contour plot of $\eta$}
        \end{subfigure}
        \begin{subfigure}{0.49\textwidth}
            \centering
            \includegraphics[width=0.95\textwidth,trim={1.2cm 6cm 1.2cm 6cm},clip]{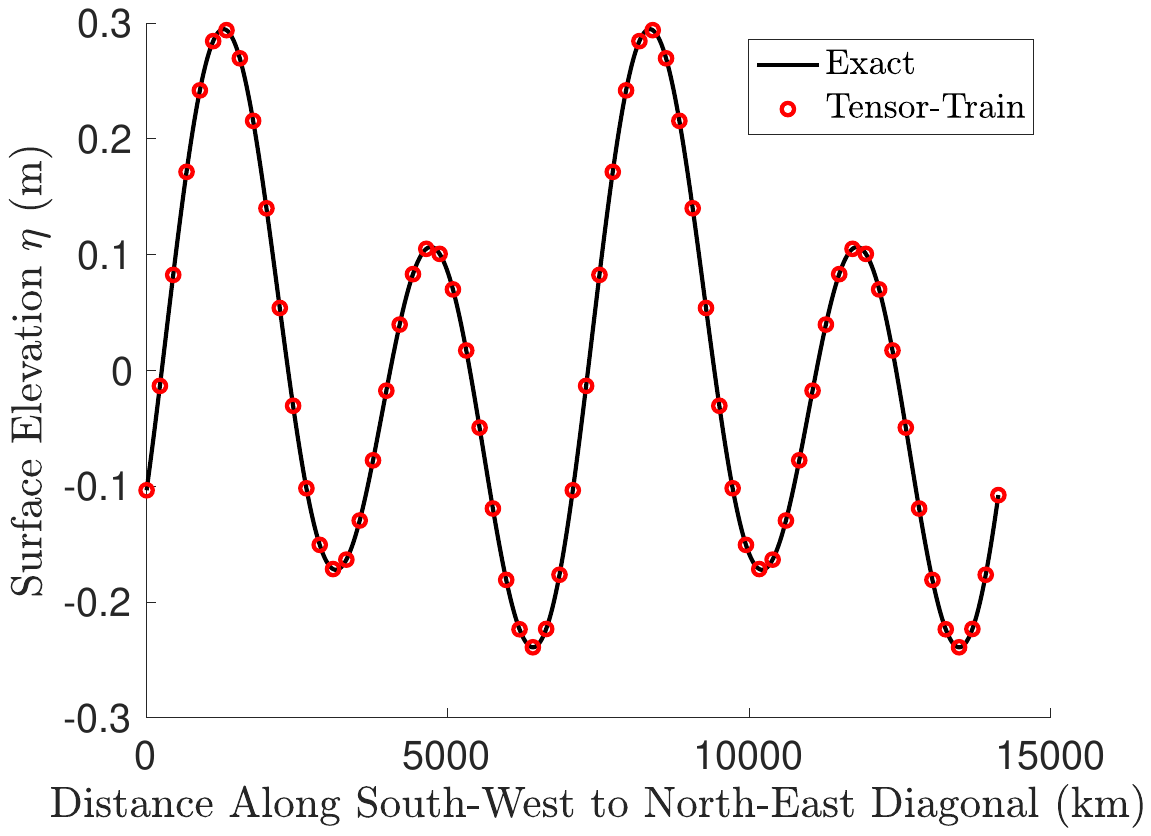}
            \caption{$\eta$ profile}
        \end{subfigure}
        \caption{Exact and TT solutions of surface elevation for Inertia-Gravity Wave at $T=3$ hours. The TT solution is obtained with Upwind5 on $1280\times1280$ grid.}\label{fig:inertia-gravity-profiles}
    \end{center}
\end{figure}

Inertia-gravity waves are an important component of geophysical turbulence \cite{Young2021}. They occur in the open ocean, and are gravity waves where particles oscillate in an ellipse due to the rotation of the Earth. The idealized test problem is linear with a flat bottom, like the Coastal Kelvin Wave, but the domain is periodic in both directions.
Specifically, we solve the linear SWEs \Cref{eq:SWE-system,eq:linear-SWE} up to the final time $T=3$ hours with $L=10^7\;m$. We consider the general solution,
\begin{equation}
    \begin{aligned}
        \eta&=\hat{\eta}\cos{\left(k_xx+k_yy-\hat{\omega}t\right)},\\
        u &= \frac{g\hat{\eta}}{\hat{\omega}^2-f^2}\left\{\hat{\omega}k_x\cos{\left(k_xx+k_yy-\hat{\omega}t\right)}-fk_y\sin{\left(k_xx+k_yy-\hat{\omega}t\right)}\right\},\\
        v &= \frac{g\hat{\eta}}{\hat{\omega}^2-f^2}\left\{\hat{\omega}k_y\cos{\left(k_xx+k_yy-\hat{\omega}t\right)}+fk_x\sin{\left(k_xx+k_yy-\hat{\omega}t\right)}\right\},
    \end{aligned}
\end{equation}
where $\hat{\omega}=\sqrt{c^2\left(k_x^2+k_y^2\right)+f^2}$. The particular solution plotted in \Cref{fig:inertia-gravity-profiles} is constructed to be the superposition of two wave modes for $\eta^{(1)}=\eta^{(2)}/2=10^{-1}$ m, $k^{(1)}_x=k^{(2)}_x/2=2\pi/L$, and $k^{(1)}_y=k^{(2)}_y/2=2\pi/L$ and initial conditions are set from these at $t=0$. 

As in the previous example, \Cref{eq:eps-tt-formula,eq:modified-eps-tt-formula} uses $C_\varepsilon=1$ for both upwind methods, but it sets $C_\varepsilon=500$ for the WENO5 scheme in this example. Time step sizes are also set as $\Delta t/\Delta x=10^{-4}$ for the Upwind3 method and $\Delta t/\Delta x^{5/3}=10^{-3}$ for the $5^{th}$ order methods. 

\begin{figure}[htbp]
    \begin{center}
        \begin{subfigure}{0.49\textwidth}
            \centering
            \includegraphics[width=1\textwidth,trim={1cm 6cm 1cm 6cm},clip]{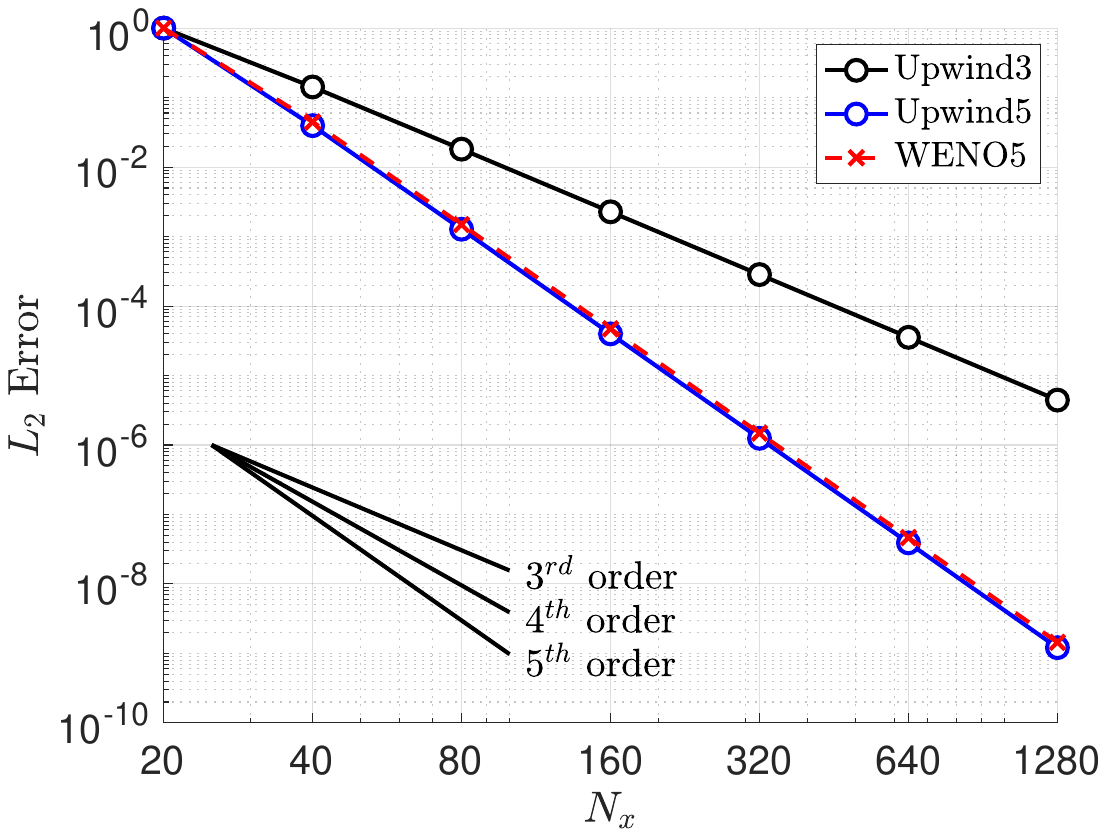}
            \caption{$L_2$ error for $\eta$}
        \end{subfigure}
        \begin{subfigure}{0.49\textwidth}
            \centering
            \includegraphics[width=0.95\textwidth,trim={1.2cm 6cm 1.2cm 6cm},clip]{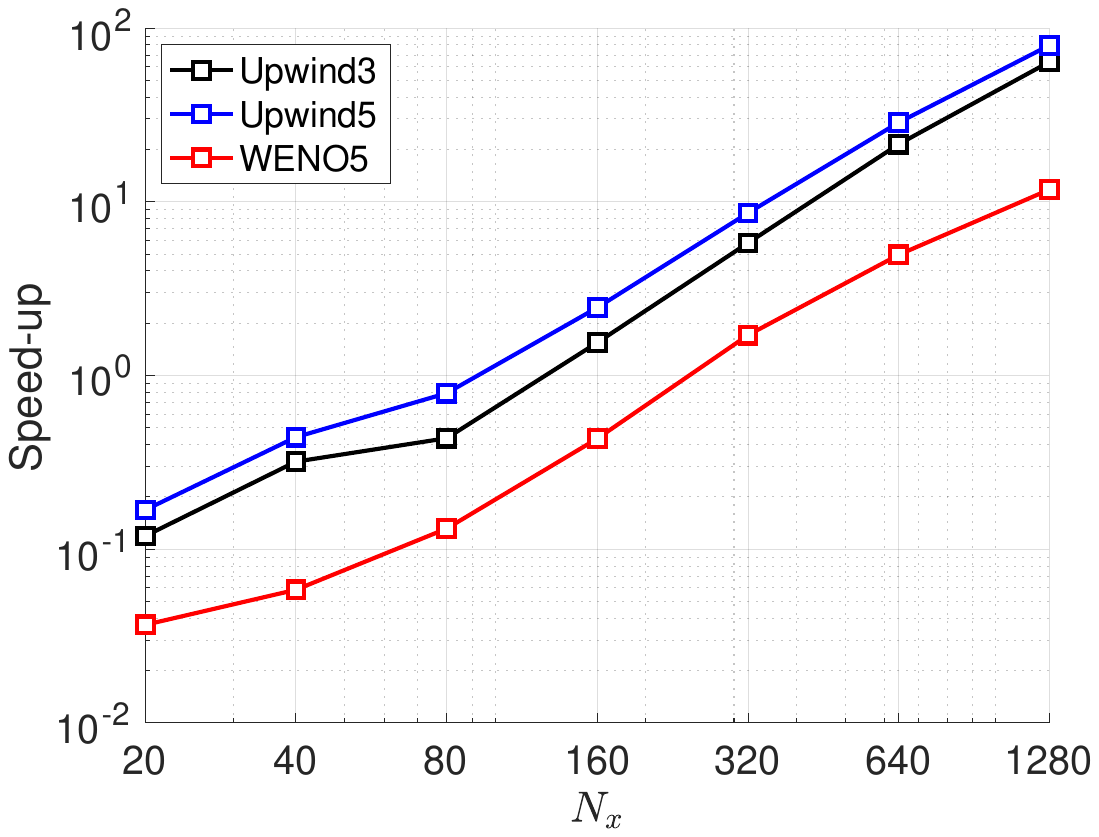}
            \caption{Speed-up}
        \end{subfigure}
        \caption{Speed-up and $L_2$ errors of surface elevation for Inertia-Gravity Wave at $T=3$ hours}\label{fig:inertia-gravity-performance}
    \end{center}
\end{figure}

In \Cref{fig:inertia-gravity-performance}, the $L_2-$convergence of $\eta$ and speed-up of the TT solvers are shown. All TT solvers recover their formal order of accuracy as in the previous example. Similar to the previous test case, TT-Upwind5 results in a slightly better speed-up than TT-Upwind3, and TT-WENO5 gives considerably less acceleration. At the finest grid level, the speed-up achieved by the TT solvers are 79x for Upwind5, 64x for Upwind3, and 12x for WENO5.
\subsection{Barotropic Tide}
\begin{figure}[htbp]
    \begin{center}
        \begin{subfigure}{0.49\textwidth}
            \centering
            \includegraphics[width=1\textwidth,trim={1cm 6cm 1cm 6cm},clip]{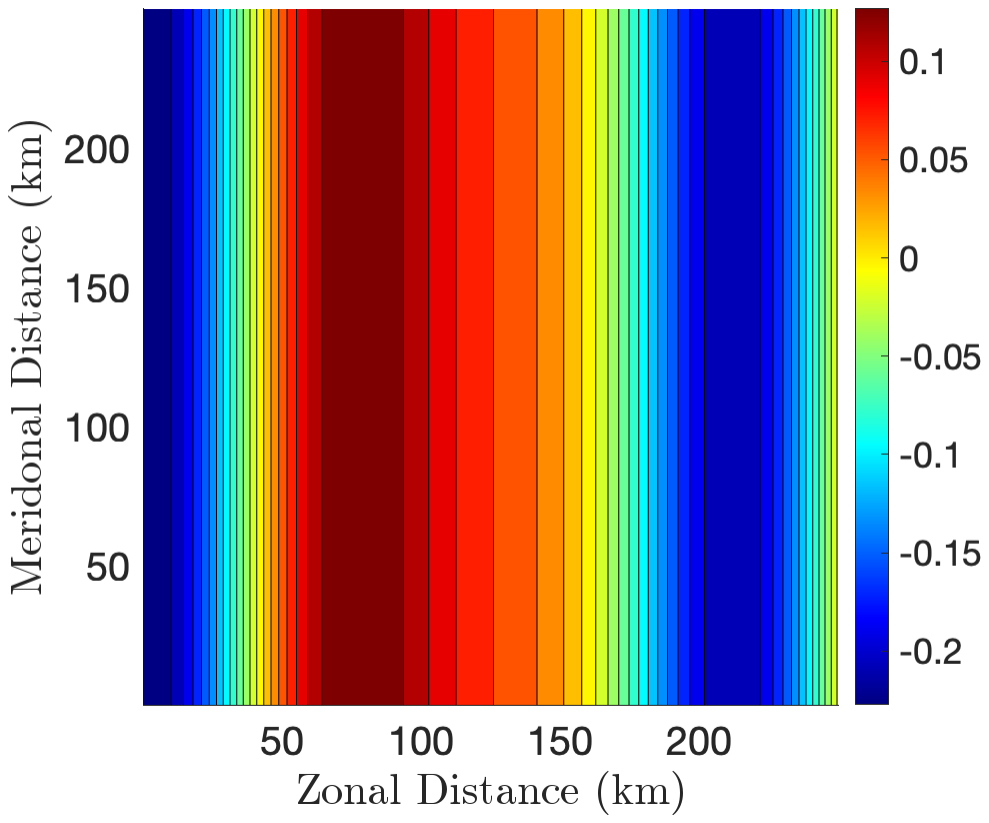}
            \caption{Contour plot of $\eta$}
        \end{subfigure}
        \begin{subfigure}{0.49\textwidth}
            \centering
            \includegraphics[width=0.95\textwidth,trim={1.2cm 6cm 1.2cm 6cm},clip]{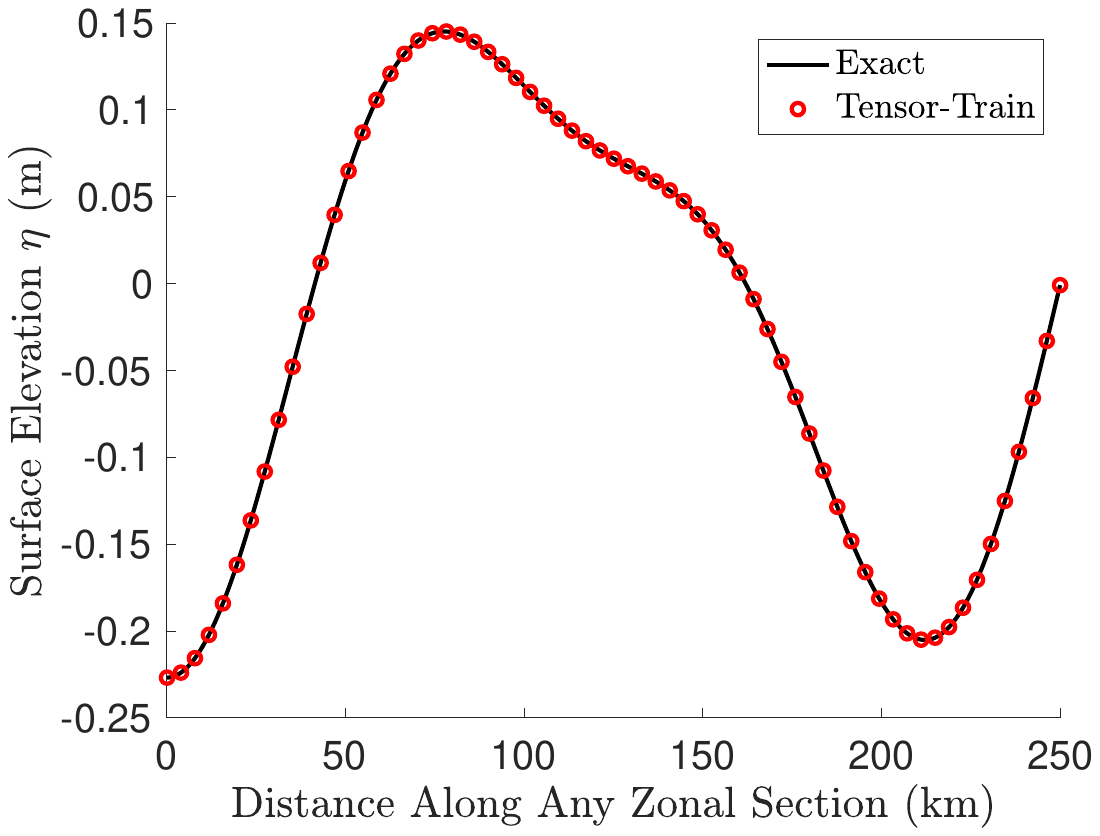}
            \caption{$\eta$ profile}
        \end{subfigure}
        \caption{Exact and TT solutions of surface elevation for Barotropic Tide at $T=30$ minutes. The TT solution is obtained with Upwind5 on $1280\times1280$ grid.}\label{fig:barotropic-tide-profiles}
    \end{center}
\end{figure}

The Barotropic Tide case is an idealized representation of a coastal tide on a continental shelf \cite{CLARKE1981665}.   Like the previous cases it tests the linear SWEs \Cref{eq:SWE-system,eq:linear-SWE}, but now on a doubly non-periodic domain.  We solve to a final time of $T=30$ minutes, where $L=25\times10^4\;m$. The general solution considered here is given as,
\begin{equation}
    \begin{aligned}
        \eta&=\hat{\eta}\cos{(kx)}\cos{(\omega t)},\\
        u &= \frac{g\hat{\eta}\omega k}{\hat{\omega}^2-f^2}\sin(kx)\sin{(\omega t)},\\
        v &= \frac{g\hat{\eta}f k}{\hat{\omega}^2-f^2}\sin(kx)\cos{(\omega t)},
    \end{aligned}
\end{equation}
where $\omega=\sqrt{gHk^2+f^2}$. We again consider a particular solution shown in \Cref{fig:barotropic-tide-profiles} constructed as a superposition of two wave modes, but now we set $\eta^{(1)}=\eta^{(2)}/2=0.2$ m, $k^{(1)}=2\pi/\lambda^{(1)}$ and $k^{(2)}=2\pi/\lambda^{(2)}$, where $\lambda^{(1)}=4L/5$ and $\lambda^{(2)}=4L/9$. Note that $H=200\;m$ in this example. Physically, the wavelengths are chosen to satisfy the conditions for tidal resonance in this domain\cite{bishnu2024}. The initial conditions at $t=0$ are set from the exact solution while the boundary condition in the $x$-direction is taken from the exact solution and is assumed periodic in the $y$-direction. 

In this test case, we set $C_\varepsilon=1$ for all reconstruction schemes, and consider $\Delta t/\Delta x = 25\times 10^{-5}$ for Upwind3 and $\Delta t/\Delta x^{5/3} = 5\times 10^{-4}$ for the $5^{th}$ order methods.

\begin{figure}[htbp]
    \begin{center}
        \begin{subfigure}{0.49\textwidth}
            \centering
            \includegraphics[width=1\textwidth,trim={1cm 6cm 1cm 6cm},clip]{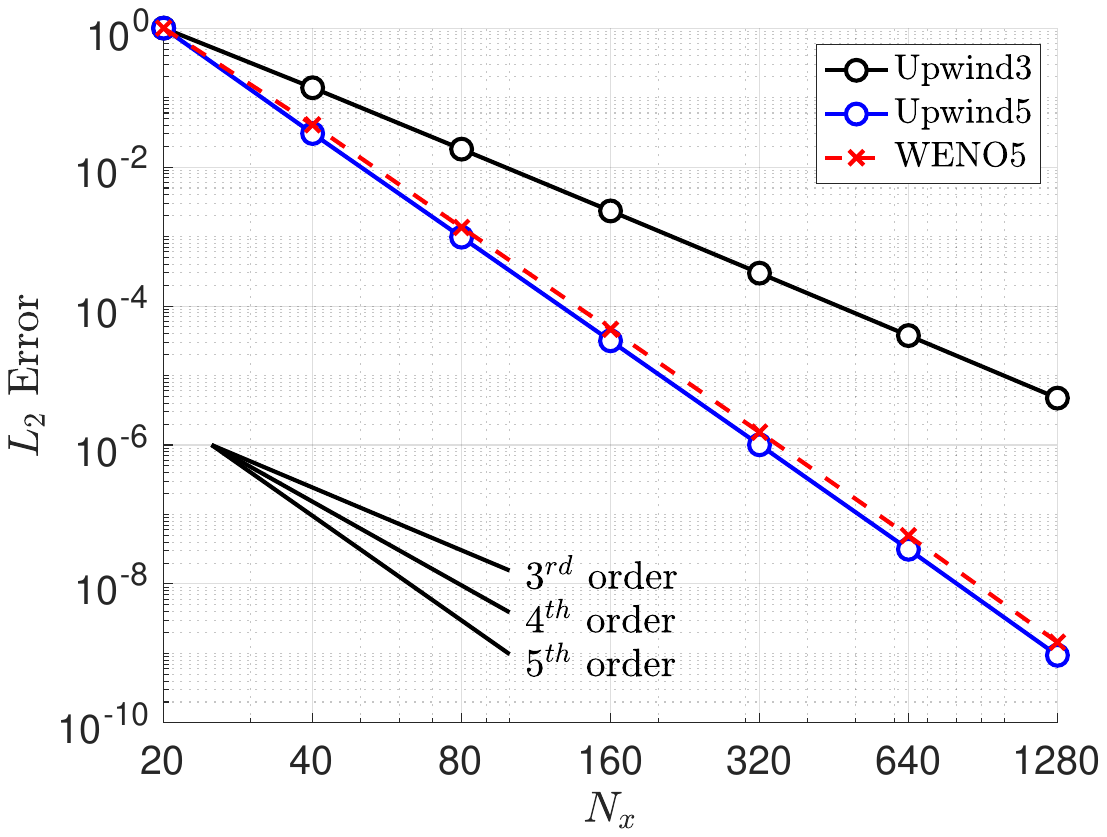}
            \caption{$L_2$ error for $\eta$}
        \end{subfigure}
        \begin{subfigure}{0.49\textwidth}
            \centering
            \includegraphics[width=0.95\textwidth,trim={1.2cm 6cm 1.2cm 6cm},clip]{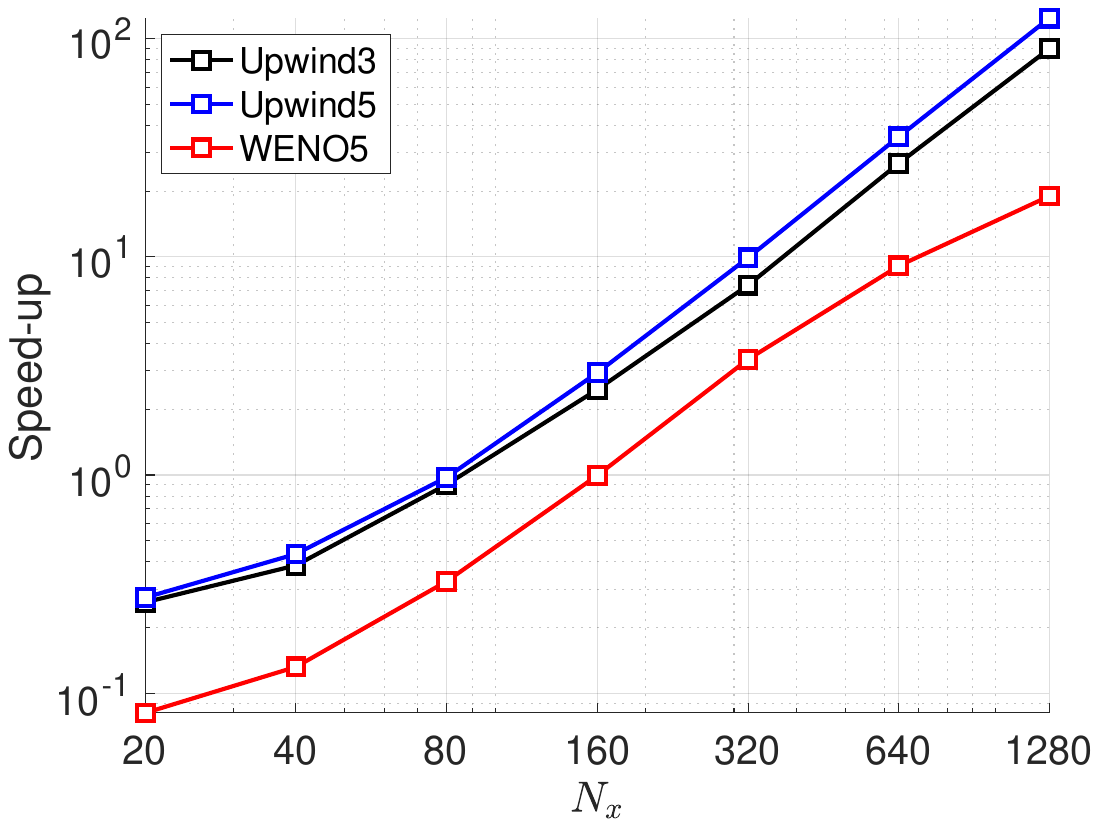}
            \caption{Speed-up}
        \end{subfigure}
        \caption{Speed-up and $L_2$ errors of surface elevation for Barotropic Tide at $T=30$ minutes}\label{fig:barotropic-tide-performance}
    \end{center}
\end{figure}

The results reported in \Cref{fig:barotropic-tide-performance} follow the trends already observed in the previous linear test cases: All TT solvers achieve their formal order of accuracy, the highest speed-up is obtained for Upwind5 and the WENO5 gives the lowest. At the finest grid level, the speed-up achieved by the TT solvers are 124x for Upwind5, 89x for Upwind3, and 19x for WENO5.

\subsection{Manufactured Solution}
\begin{figure}[htbp]
    \begin{center}
        \begin{subfigure}{0.49\textwidth}
            \centering
            \includegraphics[width=1\textwidth,trim={1cm 6cm 1cm 6cm},clip]{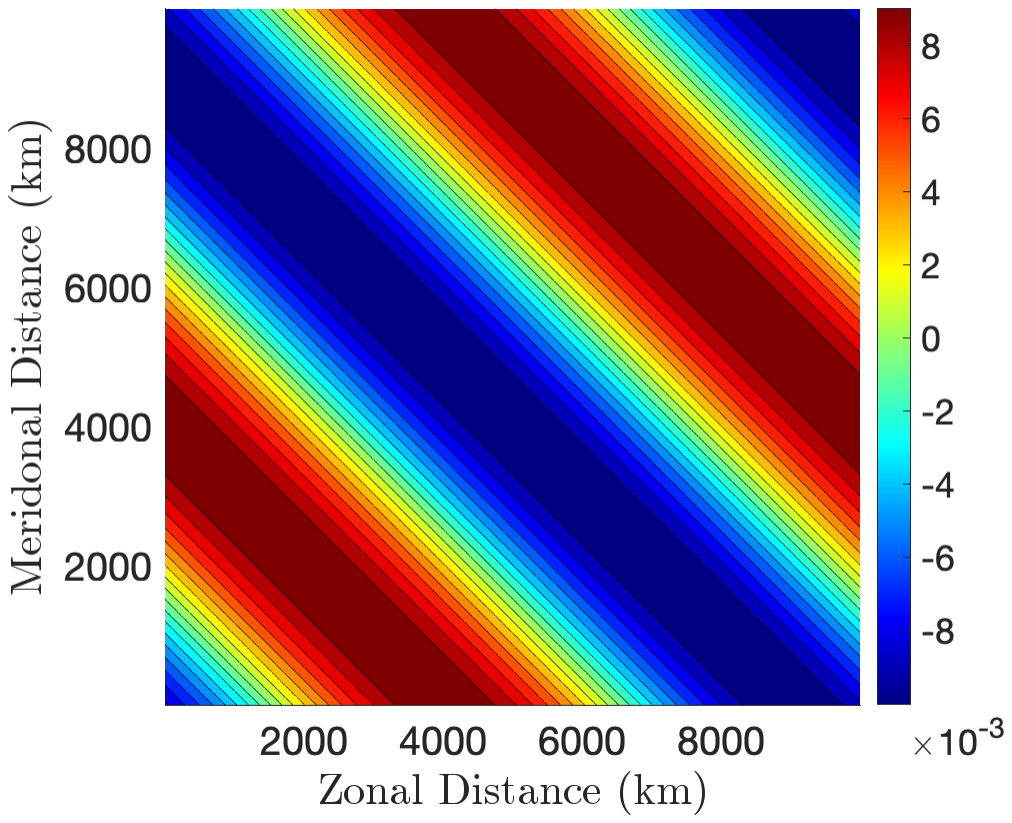}
            \caption{Contour plot of $\eta$}
        \end{subfigure}
        \begin{subfigure}{0.49\textwidth}
            \centering
            \includegraphics[width=0.95\textwidth,trim={1.2cm 6cm 1.2cm 6cm},clip]{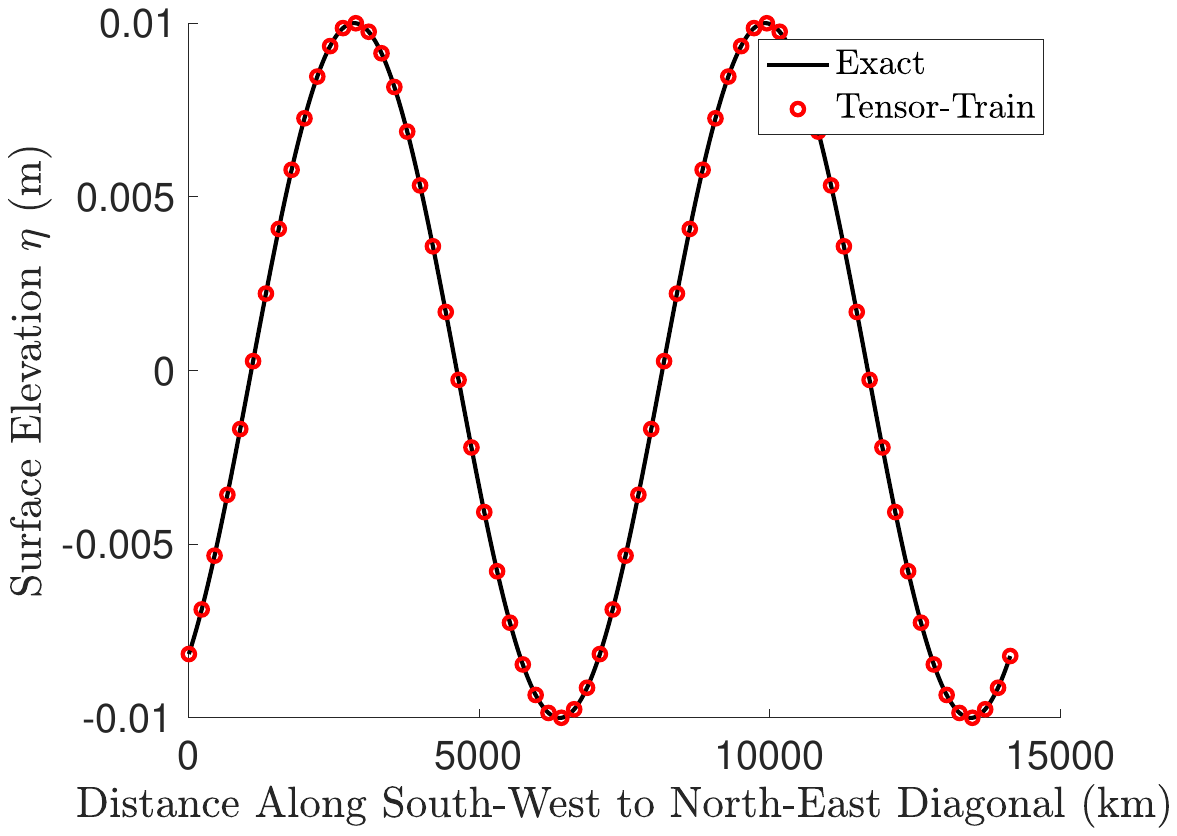}
            \caption{$\eta$ profile}
        \end{subfigure}
        \caption{Exact and TT solutions of surface elevation for the Manufactured Solution at $T=3$ hours. The TT solution is obtained with Upwind5 on $1280\times1280$ grid.}\label{fig:manuf-profiles}
    \end{center}
\end{figure}

In this test case, we turn our attention to the nonlinear SWEs and solve \Cref{eq:SWE-system,eq:nonlinear-SWE} up to the final time $T=3$ hours, where we set $L=10^7\;m$. The manufactured solution does not have an analog observed in nature, but is important because an analytic solution may be derived that tests all terms in the SWEs. Subsequent SWE test cases for geophysical turbulence with more complex behavior do not have analytic solutions\cite{williamson1992standard}. Using the method of manufactured solutions\cite{Roache2019}, we enforce
\begin{equation}
    \begin{aligned}
        \eta &= \hat{\eta}\sin{\left(k_xx+k_yy-\hat{\omega}t\right)},\\
        u &= \hat{u}\cos{\left(k_xx+k_yy-\hat{\omega}t\right)},\\
        v&=0,
    \end{aligned}
\end{equation}
where $\hat{\eta}=10^{-2}\;m$, $\hat{u}=10^{-2}\;m/s$, $k_x=k_y=2\pi/L$, and $\hat{\omega}=\sqrt{c^2\left(k_x^2+k_y^2\right)}$, as depicted in \Cref{fig:manuf-profiles}. We consider periodic boundaries and impose the initial condition from the manufactured solution at $t=0$. 
In this nonlinear test case, we set $C_\varepsilon=10^{-4}$ for the Upwind3 scheme to maintain its $3^{rd}$ order accuracy. For the $5^{th}$ order methods, however, $C_\varepsilon=1$ is deemed to be sufficient. The time stepping is performed according to $\Delta t/\Delta x = 5\times 10^{-5}$ for Upwind3 and $\Delta t/\Delta x^{5/3} = 5\times 10^{-3}$ for the $5^{th}$ order methods.

\begin{figure}[htbp]
    \begin{center}
        \begin{subfigure}{0.49\textwidth}
            \centering
            \includegraphics[width=1\textwidth,trim={1cm 6cm 1cm 6cm},clip]{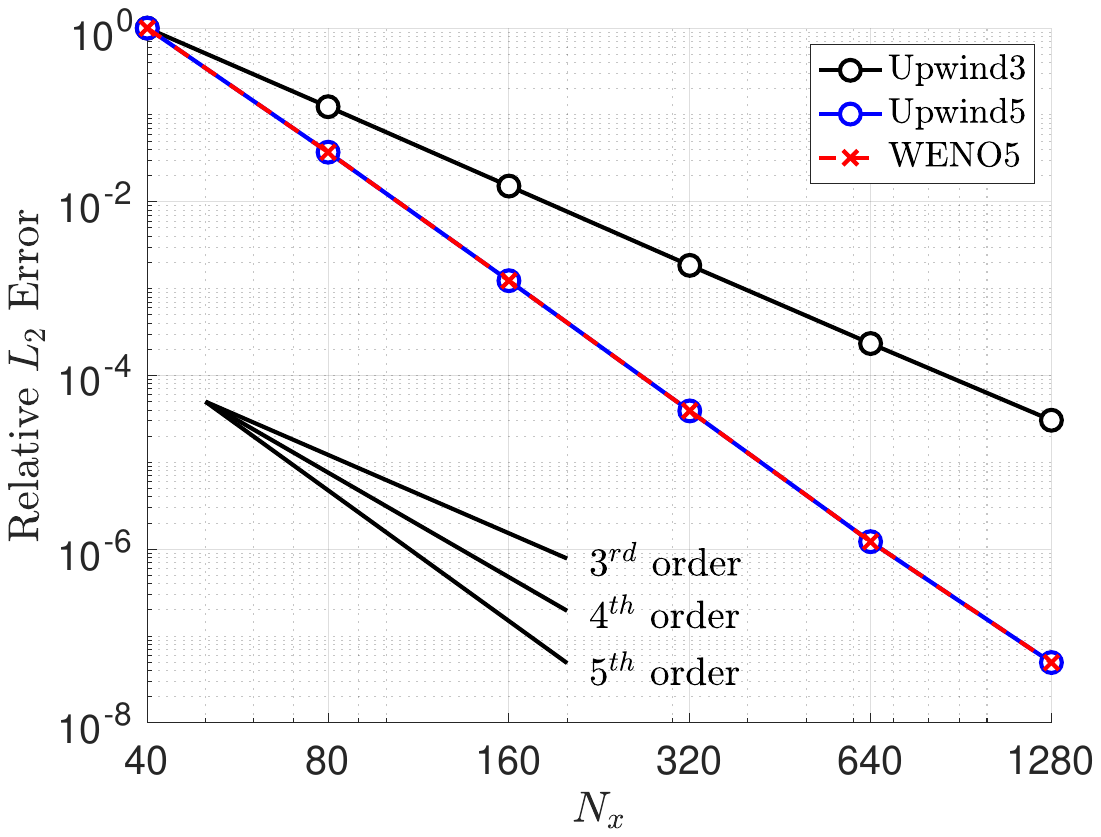}
            \caption{$L_2$ error for $h$}
        \end{subfigure}
        \begin{subfigure}{0.49\textwidth}
            \centering
            \includegraphics[width=0.95\textwidth,trim={1.2cm 6cm 1.2cm 6cm},clip]{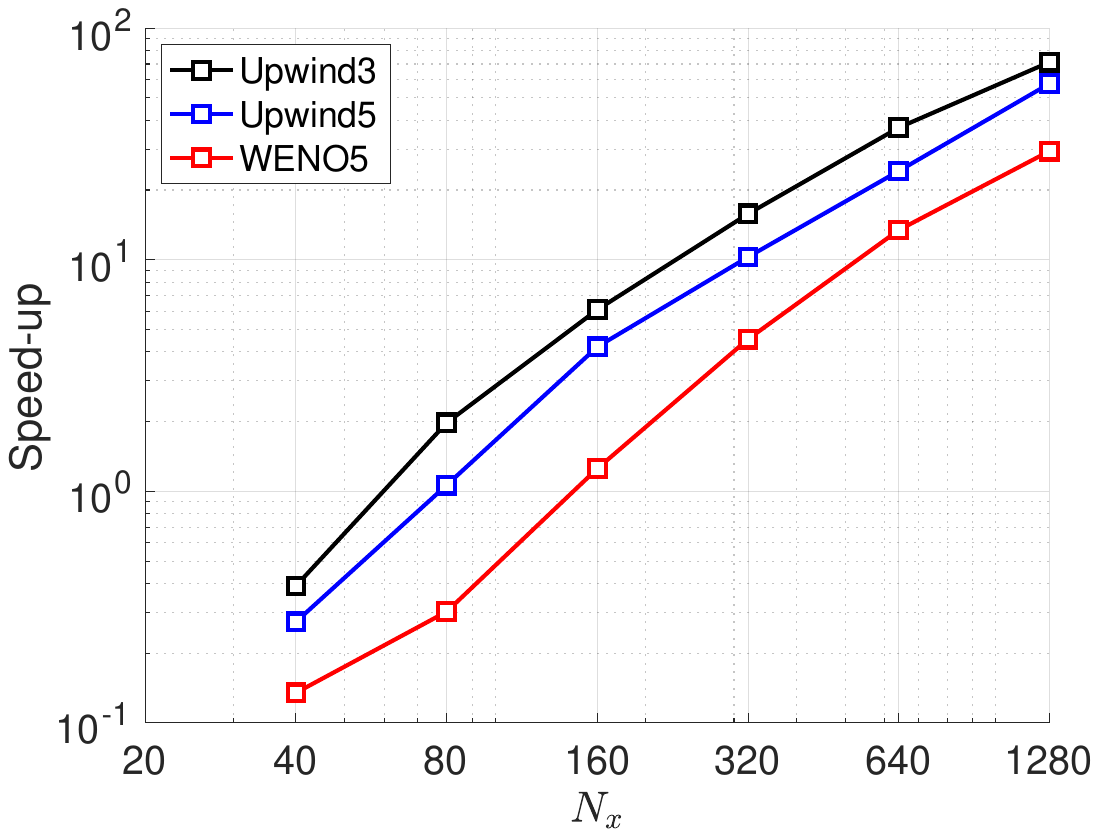}
            \caption{Speed-up}
        \end{subfigure}
        \caption{Speed-up and $L_2$ errors of the shallow water layer thickness for Manufactured Solution at $T=3$ hours.}\label{fig:manufactured-performance}
    \end{center}
\end{figure}

In \Cref{fig:manufactured-performance}, we observed the formal $3^{rd}$ order $L_2$ convergence of Upwind5 and $5^{th}$ order convergence of Upwind5 and WENO5, as in the linear cases. However, we see that Upwind3 gives a better acceleration than Upwind5 unlike the linear test cases. This may be due to the nonlinear nature of the fluxes in \Cref{eq:nonlinear-SWE}. Recall that TT rounding operation is performed after almost each TT multiplication and addition operation to avoid the rank growth. Therefore, TT rounding is applied several times for the nonlinear fluxes while no rounding is needed for the linear fluxes. Since higher-order schemes require more quadrature points to compute fluxes, the additional cost of TT rounding is naturally more expensive for higher-order method for. In addition, the nonlinear fluxes need to compute $1/h_{TT}$ using \Cref{alg:taylor-inverse}, which is also needed to be carried out on larger arrays for higher-order schemes. As a result, the different mechanics of flux computation in linear and nonlinear problems lead to different trends in terms of the speed-up. At the finest grid level, the speed-up achieved by the TT solvers are 71x for Upwind3, 57x for Upwind5, and 29x for WENO5.

%
%
\section{Conclusion}\label{sec:conclusion}
In this study, we developed a high-order tensor-train finite volume method for linear and nonlinear Shallow Water Equations (SWEs). The implementation of the $3^{rd}$ order Upwind3, $5^{th}$ order Upwind5 and order WENO5 reconstruction schemes in the tensor-train format were presented in detail. We demonstrated that all reconstruction schemes achieved their formal order of convergence in the TT format for the linear and nonlinear SWEs. The upwind methods perform reconstruction procedures that directly modify TT cores while the WENO scheme uses TT cross interpolation, due to which the WENO scheme was the method with lowest speed-up values among all numerical results. For the linear SWEs, the fluxes are computed directly without using any TT rounding or cross interpolation, which resulted in speed-up values up to 124x compared to the standard full-tensor implementation. However, the nonlinear SWEs need to calculate the reciprocal of the shallow water layer thickness (SWLT) in the TT format followed by applying TT rounding several times, which was reflected in our numerical results by reduced the speed-up values compared to the linear case. To compute the TT reciprocal of the SWLT, we employed an approximation based on Taylor series, which turned out to be quite efficient since the SWLT can usually be decomposed as a large mean value and small oscillations. Overall, we showed that the TT finite volume method maintains the accuracy of the underlying numerical discretization while significantly accelerating it for the SWEs.

The next steps will be to test this TT method with more realistic simulations---first with geophysical turbulence on the sphere with the Williamson SWE test cases \cite{williamson1992standard}, and then with primitive equations for the atmosphere and ocean. This study used Matlab for expediency, but future work will leverage high-performance computing hardware and libraries so that TT may be compared directly against standard codes and methods. TT achieves speed-up by using an SVD to truncate the numerical solution, and then time-steps a reduced version of the model state. An open question is whether TT can obtain significant speed-up with complex geophysical flows while maintaining a high-quality solution.  This study and previous studies have shown speed-ups of 100 to 1000 \cite{danis2024tensortrain} for various TT applications, which provides a strong motivation to test more complex and realistic configurations for weather and climate modeling.

\section*{Acknowledgments}
The authors gratefully acknowledge the support of the Laboratory
Directed Research and Development (LDRD) program of Los Alamos
National Laboratory under project numbers 20230067DR and 20240782ER and ISTI Rapid Response 20248215CT-IST.
Los Alamos National Laboratory is operated by Triad National Security,
LLC, for the National Nuclear Security Administration of
U.S. Department of Energy (Contract No.\ 89233218CNA000001).

%
\section*{Appendix}
In this Appendix, we will briefly describe the `nondimensionalization' of the governing equations. Denoting dimensional variables by a star, let us rewrite the Shallow Water Equations:
\begin{equation}\label{eq:dim-SWE-system}
    \frac{\partial \Q^*}{\partial t^*} + \frac{\partial \Fc^*}{\partial x^*} + \frac{\partial \Gc^*}{\partial y^*} = \Sc^*.
\end{equation}
For the linear case, we have
\begin{equation}\label{eq:dim-linear-SWE}
    \Q^* = 
    \begin{pmatrix}
        \eta^* \\
        u^*    \\
        v^*    \\
    \end{pmatrix},
    \qquad
    \Fc^* = 
    \begin{pmatrix}
        H^*u^* \\
        g^*\eta^*  \\
        0      \\
    \end{pmatrix},
    \qquad
    \Gc^* = 
    \begin{pmatrix}
        H^*v^* \\
        0      \\
        g^*\eta^*  \\
    \end{pmatrix},
    \qquad
    \Sc^* = 
    \begin{pmatrix}
        0     \\
        f^*v^*    \\
       -f^*u^*    \\
    \end{pmatrix},
\end{equation}
and for the nonlinear case,
\begin{equation}\label{eq:dim-nonlinear-SWE}
    \Q^* = 
    \begin{pmatrix}
        h^* \\
        h^*u^*    \\
        h^*v^*    \\
    \end{pmatrix},
    \qquad
    \Fc^* = 
    \begin{pmatrix}
        h^*u^* \\
        h^*{u^*}^2+\frac{1}{2}g^*{h^*}^2  \\
        h^*u^*v^*      \\
    \end{pmatrix},
    \qquad
    \Gc^* = 
    \begin{pmatrix}
        h^*v^* \\
        h^*u^*v^*      \\
        h^*{v^*}^2+\frac{1}{2}g^*{h^*}^2  \\
    \end{pmatrix},
    \qquad
    \Sc^* = 
    \begin{pmatrix}
        0   \\
        f^*h^*v^* \\
       -f^*h^*u^* \\
    \end{pmatrix}.
\end{equation}

Next, we define characteristic scales for each variable in the SWEs, such as the reference length $L^*_{ref}$, reference velocity $U^*_{ref}$, reference time $t^*_{ref}=L^*_{ref}/U^*_{ref}$, reference height $H^*_{ref}$, reference acceleration $g^*_{ref}=U^{*^2}_{ref}/H^*_{ref}$, reference frequency $f^*_{ref}=U^*_{ref}/L^*_{ref}$. Then, the nondimensional parameters are defined as
\begin{equation}
    \begin{aligned}
        &t=\frac{t^*}{t^*_{ref}},\quad x=\frac{x^*}{L^*_{ref}},\quad y=\frac{y^*}{L^*_{ref}},\\
        &\eta = \frac{\eta^*}{H^*_{ref}},\quad h = \frac{h^*}{H^*_{ref}},\\
        &u = \frac{u^*}{U^*_{ref}},\quad v = \frac{v^*}{U^*_{ref}}, \\
        &g = \frac{g^*}{g^*_{ref}},\quad f = \frac{f^*}{f^*_{ref}}\\
    \end{aligned}
\end{equation}
Substituting these into \Cref{eq:dim-SWE-system,eq:dim-linear-SWE,eq:dim-nonlinear-SWE} gives \Cref{eq:SWE-system,eq:linear-SWE,eq:nonlinear-SWE}. Note that both set of equations have the same form. However, in the dimensional form, $C_\varepsilon$ in \Cref{eq:eps-tt-formula,eq:modified-eps-tt-formula} needs to attain extremely small values, e.g. $C_\varepsilon=10^{-18}-10^{-22}$, due to extremely large computational domains. On the other hand, the nondimensional form significantly decreases this ambiguity and $C_\varepsilon=1$ seems to work well in the most numerical examples discussed in \Cref{sec:results}.  

Finally, we summarize the reference values used in the nondimensionalization of the governing equations in \Cref{tab:reference-values} below.
\begin{table}[ht]
    \centering
    \begin{tabular}{lrrr}\hlineB{3}
        Test Case             & $L^*_{ref}\;(m)$ & $U^*_{ref}\;(m/s)$     & $H^*_{ref}\;(m)$ \\\hlineB{3}
        Coastal Kelvin Wave   & $5\times10^6$  &        $5\times10^{-3}$           &       $0.1$     \\
        Inertia-Gravity Wave  &    $10^7$      & $1.622\times10^{-3}$ &       $0.2$    \\
        Barotropic Tide       & $25\times10^4$ & $3.163\times10^{-3}$ &       $0.2$    \\
        Manufactured Solution &    $10^7$      &       $10^{-2}$         &       $500$    \\\hlineB{3}
    \end{tabular}
    \caption{Characteristic Scales for Numerical Examples}
    \label{tab:reference-values}
\end{table}